%% file: 2011-03_FGM_vib.tex
\documentclass[preprint,12pt]{elsarticle}

\usepackage{amssymb,amsmath}
\usepackage{graphicx}
\usepackage{graphics}
\usepackage{subfigure}
\usepackage{mathrsfs}
\usepackage{eucal}
\usepackage{amsmath}
\usepackage{amssymb}
\usepackage{amsfonts}
\usepackage{multirow}
\usepackage{lscape}

\newcommand{\Eref}[1]{Equation (\ref{#1})}
\newcommand{\fref}[1]{Figure (\ref{#1})}
\newcommand{\Erefs}[1]{Equations (\ref{#1})}
\newcommand{\frefs}[1]{Figures~(\ref{#1})}

\newcommand{\xx}{\mathbf{x}}
\newcommand{\KK}{\mathbf{K}}
\newcommand{\bm}{\mathbf{M}}
\newcommand{\bn}{\mathbf{N}}
\newcommand{\DD}{\mathbf{D_b}}

\newcommand{\bveps}{\boldsymbol{\varepsilon}}
\newcommand{\uu}{\mathbf{u}}
\newcommand{\BB}{\mathbf{B}}
\newcommand{\qq}{\mathbf{q}}

\journal{xxxxx}

\begin{document}

\begin{frontmatter}

\title{Natural frequencies of cracked functionally graded material plates by the extended finite element method}

\author[a]{S~Natarajan}
\author[b]{P M Baiz}
\author[c]{S~Bordas\fnref{label2}\corref{cor1}}
\author[d]{T~Rabczuk}
\author[e]{P~Kerfriden}

\address[a]{PhD Research Student, Institute of Mechanics and Advanced Materials, Theoretical and Computational Mechanics,
Cardiff University, U.K.}
\address[b]{Lecturer, Department of Aeronautics, Imperial College, London, U.K.}
\address[c]
{Professor, Institute of Mechanics and Advanced Materials, Theoretical and Computational Mechanics,
Cardiff University, U.K.}
\address[d]
{Professor, Department of Civil Engineering, Bauhaus-Universit\"{a}t
Weimar, Germany}
\address[e]
{Lecturer, Institute of Mechanics and Advanced Materials, Theoretical and Computational Mechanics,
Cardiff University, U.K.}

\fntext[4]{Cardiff School of Engineering, Theoretical, Applied and Computational Mechanics, Cardiff University, Wales, U.K. Email: stephane.bordas@alumni.northwestern.edu.
Tel. +44 (0)29 20875941.}

\begin{abstract}
In this paper, the linear free flexural vibration of cracked functionally graded material plates is studied using the extended finite element method. A 4-noded quadrilateral plate bending element based on field and edge consistency requirement with 20 degrees of freedom per element is used for this study. The natural frequencies and mode shapes of simply supported and clamped square and rectangular plates are computed as a function of gradient index, crack length, crack orientation and crack location. The effect of thickness and influence of multiple cracks is also studied.
\end{abstract}


\end{frontmatter}

\section{Introduction}
Engineered materials such as laminated composites are widely used in automotive and aerospace industry due to their excellent strength-to and stiffness-to-weight ratios and their possibility of tailoring their properties in optimizing their structural response. But due to sudden change in material properties between the layers in laminated composites, these materials suffer from premature failure or by the decay of stiffness characteristics because of delaminations and chemically unstable matrix and lamina adhesives. The emergence of functionally graded materials (FGMs)~\cite{Koizumi1993,Koizumi1997} has revolutionized the aerospace and aerocraft industry. The FGMs used initially as thermal barrier materials for aerospace structural applications and fusion reactors are now developed for general use as structural components in high temperature environments. FGMs are manufactured by combining metals and ceramics. These materials are inhomogeneous, in the sense that the material properties vary smoothly and continuously in one or more directions. FGMs combine the best properties of metals and ceramics and are strongly considered as a potential structural material candidates for certain class of aerospace structures exposed to a high temperature environment.

It is seen from the literature that the amount of work carried out on the vibration characteristics of FGMs is considerable~\cite{Vel2004,Matsunaga2008,He2001,Liew2001,Ng2000,Yang2002,Ferreira2006}. He \textit{et al.,}~\cite{He2001} presented finite element formulation based on thin plate theory for the vibration control of FGM plate with integrated piezoelectric sensors and actuators under mechanical load whereas Liew \textit{et al.,}~\cite{Liew2001}have analyzed the active vibration control of plate subjected to a thermal gradient using shear deformation theory. Ng \textit{et al.,}~\cite{Ng2000} have investigated the parametric resonance of plates based on Hamilton`s principle and the assumed mode technique. Yang and Shen~\cite{Yang2001} have analyzed the dynamic response of thin FGM plates subjected to impulsive loads using Galerkin procedure coupled with modal superposition method whereas, by neglecting the heat conduction effect, such plates and panels in thermal environments have been examined based on shear deformation with temperature dependent material properties in~\cite{Yang2002}.
Qian \textit{et al.,}~\cite{Qian2004} studied the static deformation and vibration of FGM plates based on higher-order shear deformation theory using meshless local Petrov-Galerkin method. Matsunaga~\cite{Matsunaga2008} presented analytical solutions for simply supported rectangular FGM plates based on second-order shear deformation plates. Vel and Batra~\cite{Vel2004} proposed three-dimensional solutions for vibrations of simply supported rectangular plates. Reddy~\cite{Reddy2000} presented a finite element solution for the dynamic analysis of a FGM plate and Ferreira~\textit{et al.,}~\cite{Ferreira2006} peformed dynamic analysis of FGM plate based on higher order shear and normal deformable plate theory using the meshless local Petrov-Galerkin method. Akbari \textit{et al.,}~\cite{R2010} studied two-dimensional wave propagation in functionally graded solids using the meshless local Petrov-Galerkin method. The above list is no way comprehensive and interested readers are referred to the literature.

FGM plates or in general plate structures, may develop flaws during manufacturing or after they have been subjected to large cyclic loading. Hence it is important to understand the dynamic response of a FGM plate with an internal flaw. It is known that cracks or local defects affect the dynamic response of a structural member. This is because, the presence of the crack introduces local flexibility and anisotropy. Moreover the crack will open and close depending on the vibration amplitude. The vibration of cracked plates was studied as early as 1969 by Lynn and Kumbasar~\cite{Lynn1967} who used a Green's function approach. Later, in 1972, Stahl and Keer~\cite{Stahl1972} studied the vibration of cracked rectangular plates using elasticity methods. The other numerical  methods that are used to study the dynamic response and instability of plates with cracks or local defects are: (1) Finite fourier series transform~\cite{Solecki1985}; (2) Rayleigh-Ritz Method~\cite{Khadem2000}; (3) harmonic balance method~\cite{Wu2005}; and (4) finite element method~\cite{Qian1991,Lee1993}. Recently, ~\cite{Huang2011} proposed solutions for the virbations of side-cracked FGM thick plates based on Reddy third-order shear deformation theory using Ritz technique. Kitipornchai~\textit{et al.,}~\cite{Kitipornchai2009} studied nonlinear vibration of edge cracked functionally graded Timoshenko beams using Ritz method. Yang~\textit{et al.,}~\cite{Yang2010} studied the nonlinear dynamic response of a functionally graded plate with a through-width crack based on Reddy's third-order shear deformation theory using a Galerkin method.

In this paper, we apply the extended finite element method (XFEM) to study the free flexural vibrations of cracked FGM plates based on first order shear deformation theory. We carry out a parametric study on the influence of gradient index, crack location, crack length, crack orientation and thickness on the natural frequencies of FGM plates using the 4-noded shear flexible element based on field and edge consistency approach~\cite{Somashekar1987}. The effect of boundary conditions and multiple cracks is also studied. Earlier, the XFEM has been applied to study the vibration of cracked isotropic plates~\cite{Bachene2009,Bachene2009a,Tiberkak2009}. Their study focussed on center and edge cracks with simply supported and clamped boundary conditions. 

The paper is organized as follows, the next section will give an introduction to FGM and a brief overview of Reissner-Mindlin plate theory. Section~\ref{XFEM:basics} illustrates the basic idea of XFEM as applicable to plates. Section~\ref{numerics} presents results for the free flexural vibration of cracked plates under different boundary conditions and aspect ratios, followed by concluding remarks in the last section.

\section{Theoretical Formulation}
\subsection{Functionally Graded Material}
A functionally graded material (FGM) rectangular plate (length $a$, width $b$ and thickness $h$), made by mixing two distinct material phases: a metal and ceramic is considered with coordinates $x,y$ along the in-plane directions and $z$ along the thickness direction (see \fref{fig:platefig}). The material on the top surface $(z=h/2)$ of the plate is ceramic and is graded to metal at the bottom surface of the plate $(z=-h/2)$ by a power law distribution. The homogenized material properties are computed using the Mori-Tanaka Scheme~\cite{Mori1973,Benvensite1987}. 

\subsubsection*{Estimation of mechanical and thermal properties}
Based on the Mori-Tanaka homogenization method, the effective bulk modulus $K$ and shear modulus $G$ of the FGM are evaluated as~\cite{Mori1973,Benvensite1987,Cheng2000,Qian2004}

\begin{eqnarray}
{K - K_m \over K_c - K_m} &=& {V_c \over 1+(1-V_c){3(K_c - K_m) \over 3K_m + 4G_m}} \nonumber \\
{G - G_m \over G_c - G_m} &=& {V_c \over 1+(1-V_c){(G_c - G_m) \over G_m + f_1}}
\label{eqn:bulkshearmodulus}
\end{eqnarray}

where

\begin{equation}
f_1 = {G_m (9K_m + 8G_m) \over 6(K_m + 2G_m)}
\end{equation}

\begin{figure}
\includegraphics{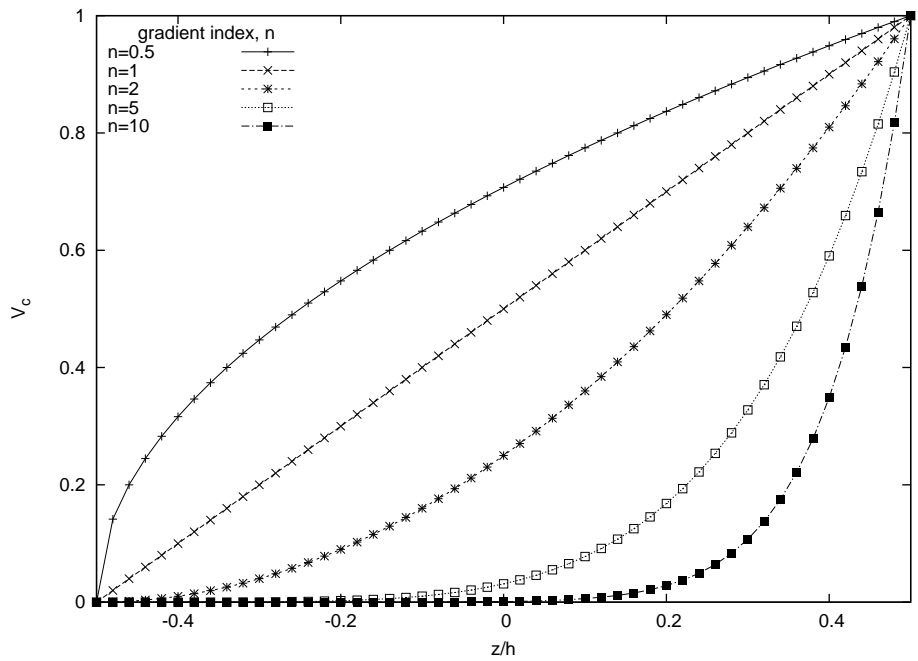}
\caption{Through thickness variation of volume fraction}
\label{fig:volfrac}
\end{figure}

Here, $V_i~(i=c,m)$ is the volume fraction of the phase material. The subscripts $c$ and $m$ refer to the ceramic and metal phases, respectively. The volume fractions of the ceramic and metal phases are related by $V_c + V_m = 1$, and $V_c$ is expressed as

\begin{equation}
V_c(z) = \left( {2z + h \over 2h} \right)^n, \hspace{0.2cm}  n \ge 0
\label{eqn:volFrac}
\end{equation}

where $n$ in \Eref{eqn:volFrac} is the volume fraction exponent, also referred to as the gradient index. \fref{fig:volfrac} shows the variation of the volume fractions of ceramic and metal, respectively, in the thickness direction $z$ for the FGM plate. The top surface is ceramic rich and the bottom surface is metal rich. 

The effective Young's modulus $E$ and Poisson's ratio $\nu$ can be computed from the following expressions:

\begin{eqnarray}
E = {9KG \over 3K+G} \nonumber \\
\nu = {3K - 2G \over 2(3K+G)}
\label{eqn:young}
\end{eqnarray}

The effective mass density $\rho$ is given by the rule of mixtures as~\cite{Vel2004}

\begin{equation}
\rho = \rho_c V_c + \rho_m V_m
\label{eqn:mdensity}
\end{equation}

The material properties
$P$ that are temperature dependent can be written as~\cite{Reddy1998}

\begin{equation}
P = P_o(P_{-1}T^{-1} + 1 + P_1 T + P_2 T^2 + P_3 T^3),
\end{equation}

where $P_o,P_{-1},P_1,P_2,P_3$ are the coefficients of temperature $T$ and are unique to each constituent material phase. 


\subsection{Plate formulation}
\begin{figure}
\centering
\input{plate.pstex_t}
\caption{Co-ordinate system of a rectangular FGM plate.}
\label{fig:platefig}
\end{figure}
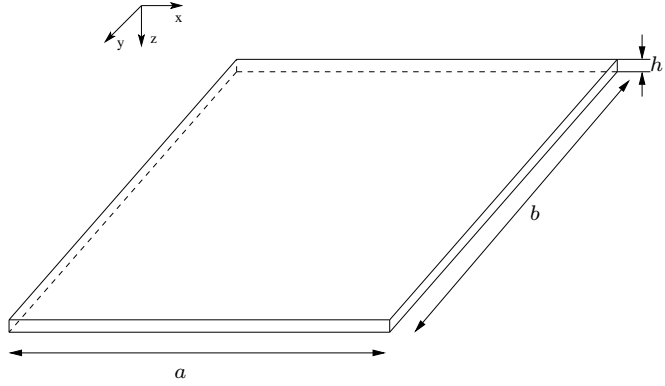

Using the Mindlin formulation, the displacements $u,v,w$ at a point $(x,y,z)$ in the plate (see \fref{fig:platefig}) from the medium surface are expressed as functions of the mid-plane displacements $u_o,v_o,w_o$ and independent rotations $\theta_x,\theta_y$ of the normal in $yz$ and $xz$ planes, respectively, as

\begin{eqnarray}
u(x,y,z,t) &=& u_o(x,y,t) + z \theta_x(x,y,t) \nonumber \\
v(x,y,z,t) &=& v_o(x,y,t) + z \theta_y(x,y,t) \nonumber \\
w(x,y,z,t) &=& w_o(x,y,t) 
\label{eqn:displacements}
\end{eqnarray}

where $t$ is the time. The strains in terms of mid-plane deformation can be written as

\begin{equation}
\bveps  = \left\{ \begin{array}{c} \bveps_p \\ 0 \end{array} \right \}  + \left\{ \begin{array}{c} z \bveps_b \\ \bveps_s \end{array} \right\} 
\label{eqn:strain1}
\end{equation}

The midplane strains $\bveps_p$, bending strain $\bveps_b$, shear strain $\bveps_s$ in \Eref{eqn:strain1} are written as

\begin{eqnarray}
\renewcommand{\arraystretch}{1.5}
\bveps_p = \left\{ \begin{array}{c} u_{o,x} \\ v_{o,y} \\ u_{o,y}+v_{o,x} \end{array} \right\}, \hspace{1cm}
\renewcommand{\arraystretch}{1.5}
\bveps_b = \left\{ \begin{array}{c} \theta_{x,x} \\ \theta_{y,y} \\ \theta_{x,y}+\theta_{y,x} \end{array} \right\} \nonumber \\
\renewcommand{\arraystretch}{1.5}
\bveps_s = \left\{ \begin{array}{c} \theta _x + w_{o,x} \\ \theta _y + w_{o,y} \end{array} \right\}, \hspace{1cm}
\renewcommand{\arraystretch}{1.5}
\end{eqnarray}

where the subscript `comma' represents the partial derivative with respect to the spatial coordinate succeeding it. The membrane stress resultants $\bn$ and the bending stress resultants $\bm$ can be related to the membrane strains, $\bveps_p$ and bending strains $\bveps_b$ through the following constitutive relations

\begin{eqnarray}
\bn &=& \left\{ \begin{array}{c} N_{xx} \\ N_{yy} \\ N_{xy} \end{array} \right\} = \mathbf{A} \bveps_p + \BB \bveps_b \nonumber \\
\bm &=& \left\{ \begin{array}{c} M_{xx} \\ M_{yy} \\ M_{xy} \end{array} \right\} = \BB \bveps_p + \DD \bveps_b  
\end{eqnarray}

where the matrices $\mathbf{A} = A_{ij}, \BB= B_{ij}$ and $\DD = D_{ij}; (i,j=1,2,6)$ are the extensional, bending-extensional coupling and bending stiffness coefficients and are defined as

\begin{equation}
\left\{ A_{ij}, ~B_{ij}, ~ D_{ij} \right\} = \int_{-h/2}^{h/2} \overline{Q}_{ij} \left\{1,~z,~z^2 \right\}~dz
\end{equation}

Similarly, the transverse shear force $Q = \{Q_{xz},Q_{yz}\}$ is related to the transverse shear strains $\bveps_s$ through the following equation

\begin{equation}
Q_{ij} = E_{ij} \bveps_s
\end{equation}

where $E_{ij} = \int_{-h/2}^{h/2} \overline{Q} \upsilon_i \upsilon_j~dz;~ (i,j=4,5)$ is the transverse shear stiffness coefficient, $\upsilon_i, \upsilon_j$ is the transverse shear coefficient for non-uniform shear strain distribution through the plate thickness. The stiffness coefficients $\overline{Q}_{ij}$ are defined as

\begin{eqnarray}
\overline{Q}_{11} = \overline{Q}_{22} = {E(z,T) \over 1-\nu^2}; \hspace{1cm} \overline{Q}_{12} = {\nu E(z,T) \over 1-\nu^2}; \hspace{1cm} \overline{Q}_{16} = \overline{Q}_{26} = 0 \nonumber \\
\overline{Q}_{44} = \overline{Q}_{55} = \overline{Q}_{66} = {E(z,T) \over 2(1+\nu) }
\end{eqnarray}

where the modulus of elasticity $E(z,T)$ and Poisson's ratio $\nu$ are given by \Eref{eqn:young}. The strain energy function $U$ is given by

\begin{equation}
U(\boldsymbol{\delta}) = {1 \over 2} \int_{\Omega} \left\{ \bveps_p^{\textup{T}} \mathbf{A} \bveps_p + \bveps_p^{\rm T} \mathbf{B} \bveps_b + 
\bveps_b^{\textup{T}} \mathbf{B} \bveps_p + \bveps_b^{\textup{T}} \mathbf{D} \bveps_b +  \bveps_s^{\textup{T}} \mathbf{E} \bveps_s\right\}~ d\Omega
\label{eqn:potential}
\end{equation}

where $\boldsymbol{\delta} = \{u,v,w,\theta_x,\theta_y\}$ is the vector of the degree of freedom associated to the displacement field in a finite element discretization. Following the procedure given in~\cite{Rajasekaran1973}, the strain energy function $U$ given in~\Eref{eqn:potential} can be rewritten as

\begin{equation}
U(\boldsymbol{\delta}) = {1 \over 2}  \boldsymbol{\delta}^{\textup{T}} \mathbf{K}  \boldsymbol{\delta}
\label{eqn:poten}
\end{equation}

where $\KK$ is the linear stiffness matrix. The kinetic energy of the plate is given by

\begin{equation}
T(\boldsymbol{\delta}) = {1 \over 2} \int_{\Omega} \left\{I_o (\dot{u}_o^2 + \dot{v}_o^2 + \dot{w}_o^2) + I_1(\dot{\theta}_x^2 + \dot{\theta}_y^2) \right\}~d\Omega
\label{eqn:kinetic}
\end{equation}

where $I_o = \int_{-h/2}^{h/2} \rho(z)~dz, ~ I_1 = \int_{-h/2}^{h/2} z^2 \rho(z)~dz$ and $\rho(z)$ is the mass density that varies through the thickness of the plate given by~\Eref{eqn:mdensity}. Substituting \Eref{eqn:poten} - (\ref{eqn:kinetic}) in Lagrange's equations of motion, the following governing equation is obtained

\begin{equation}
\bm \ddot{\boldsymbol{\delta}} + \KK \boldsymbol{\delta} = \mathbf{0}
\label{eqn:govereqn}
\end{equation}

where $\bm$ is the consistent mass matrix. After substituting the characteristic of the time function~\cite{Ganapathi1991} $\ddot{\boldsymbol{\delta}} = -\omega^2 \boldsymbol{\delta}$, the following algebraic equation is obtained

\begin{equation}
\left( \KK  - \omega^2 \bm\right) \boldsymbol{\delta} = \mathbf{0}
\label{eqn:finaldiscre}
\end{equation}

where $\KK$ is the stiffness matrix, $\omega$ is the natural frequency. 

\section{Overview of the extended finite element method}
\label{XFEM:basics} In this section, we give a brief overview of the
XFEM for plates. By exploiting the idea of partition of unity noted
by Babu\v{s}ka \textit{et al.,}~\cite{Babuvska1994}, Belytschko's
group~\cite{Belytschko1999} introduced the XFEM to solve linear
elastic fracture mechanics problems. The conventional expansion of
the displacement field using a polynomial basis fails to capture the
local behavior of the problem (e.g., steep stress gradients,
material discontinuity). Hence, new functions are added to the
conventional set of basis functions that contain information regarding the
localized behavior. In general, the field variables are approximated
by~\cite{Belytschko1999,Melenk1995,Duarte2007,Babuvska1997,Babuska2008,Simone2006}:

\begin{equation}
\uu^h(\xx)=\sum\limits_{i \in \mathcal{N}^{\rm{fem}}} {N_i
}(\xx)\qq_i + \textup{enrichment functions} \label{eqn:xfemappr}
\end{equation}

where $N_i(\xx)$ are standard finite element shape functions, $\qq_I$ are nodal variables associated with node $I$. In the following,
we briefly describe the standard discretization of a plate using the field consistent Q4 plate element.
The enriched field consistent Q4 element is then described. And finally, the discretized equations for
the eigenvalue problem is given. In this section, only the essential details are given,
interested readers are referred to recent review papers on XFEM~\cite{Karihaloo2003,Yazid2009}. A review on the implementation of the extended finite element is given in~\cite{Bordas2007a} and a detailed description of the state of the art on the simulation of cracks by partition of unity enriched methods is given in~\cite{Rabczuk2010}.

\subsection{Element Description}

The plate element employed here is a $\mathcal{C}^0$ continuous shear flexible field consistent element with five degrees of freedom $(u_o,v_o,w_o,\theta_x,\theta_y)$ at four nodes in an 4-noded quadrilateral (QUAD-4) element. If the interpolation functions for QUAD-4 are used directly to interpolate the five variables $(u_o,v_o,w_o,\theta_x,\theta_y)$ in deriving the shear strains and membrane strains, the element will lock and show oscillations in the shear and membrane stresses. The field consistency requires that the transverse shear strains and membrane strains must be interpolated in a consistent manner. Thus, the $\theta_x$ and $\theta_y$ terms in the expressions for shear strain $\bveps_s$ have to be consistent with the derivative of the field functions, $w_{o,x}$ and $w_{o,y}$. This is achieved by using field redistributed substitute shape functions to interpolate those specific terms, which must be consistent as described in~\cite{Somashekar1987,Ganapathi1991}. This element is free from locking and has good convergence properties. For complete description of the element, interested readers are referred to the literature~\cite{Somashekar1987,Ganapathi1991}, where the element behavior is discussed in great detail. Since the element is based on the field consistency approach, exact integration is applied for calculating various strain energy terms. 

\subsection{Enriched Q4 element}
Consider a mesh of field consistent Q4 elements and an independent crack geometry
as shown in~\fref{fig:xfem_crack}. The following enriched
approximation proposed by Dolbow \text{et al.,}~\cite{Dolbow2000}
for the plate displacements and the section rotations are used:

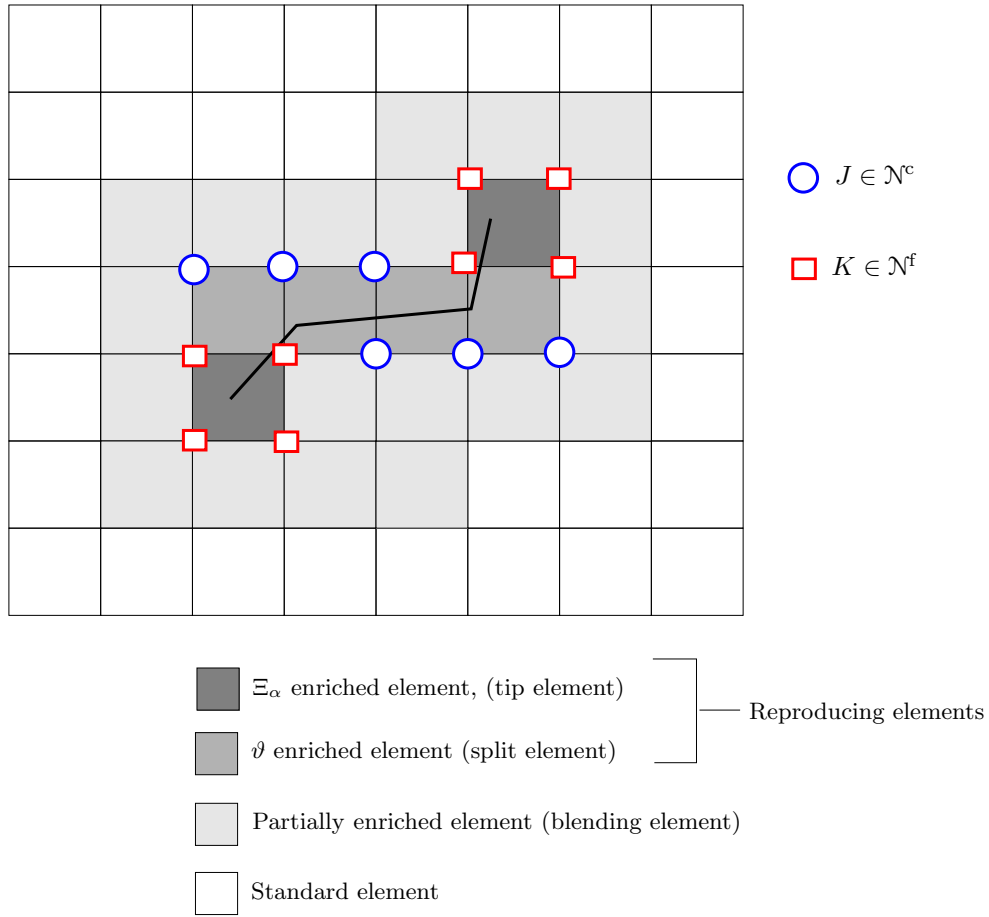
\begin{figure}[htpb]
\centering
\input{XFEMElemcate.pstex_t}
\caption{A typical FE mesh with an arbitrary crack. `Squared' nodes are enriched with the heaviside function and `circled' nodes with the near tip functions, which allows representing cracks independent of the background mesh.}
\label{fig:xfem_crack}
\end{figure}

\begin{equation}
\begin{split}
(u^h,v^h,w^h)\left(\xx\right) = \sum_{i \in \mathcal{N}^{\rm{fem}}} \phi_i(\xx) (u_i,v_i,w_i) + \sum_{j \in
\mathcal{N}^{\rm{c}}} \phi_j(\xx) H(\xx) (b_j^u,b_j^v,b_j^w) +  \\ \sum_{k \in \mathcal{N}^{\rm{f}}}\phi_k(\xx) \left(
\sum_{l=1}^4 (c_{kl}^u,c_{kl}^v,c_{kl}^w) G_l (r,\theta) \right)
\end{split}
\label{eqn:platexfem1}
\end{equation}

The section rotations in the shear terms are approximated by

\begin{equation}
\begin{split}
(\theta_x^h,\theta_y^h)\left(\xx\right) = \sum_{i \in \mathcal{N}^{\rm{fem}}} \tilde{\phi}_i(\xx) (\theta_{x_i},\theta_{y_i}) +
\sum_{j \in \mathcal{N}^{\rm{c}}} \tilde{\phi}_j(\xx) H(\xx) (b_j^{\theta_x},b_j^{\theta_y}) + \\ 
\sum_{k \in
\mathcal{N}^{\rm{f}}}\tilde{\phi}_k(\xx) \left( \sum_{l=1}^4 (c_{kl}^{\theta_x},c_{kl}^{\theta_y}) F_l
(r,\theta) \right). 
\end{split}
\label{eqn:platexfem2}
\end{equation}

In \Erefs{eqn:platexfem1} and (\ref{eqn:platexfem2}), 
$(u_i,v_i,w_i,\theta_{x_i},\theta_{y_i})$ are the nodal unknown vectors associated with the continuous part of the finite element solution, $b_j$ is the nodal enriched degree of freedom vector associated with the Heaviside (discontinuous) function, and $c_{kl}$ is the nodal enriched degree of freedom vector associated with the elastic asymptotic crack-tip functions. 
The asymptotic functions, $G_l$ and $F_l$ in
\Erefs{eqn:platexfem1} and (\ref{eqn:platexfem2}) are given by
(\cite{Dolbow2000}):

\begin{subequations}
\begin{eqnarray}
G_l(r,\theta) \equiv \left\{ \sqrt[3]{r}
\sin\left(\frac{\theta}{2}\right), \sqrt[3]{r}
\cos\left(\frac{\theta}{2}\right), \sqrt[3]{r}
\sin\left(\frac{3\theta}{2}\right), \sqrt[3]{r}
\cos\left(\frac{3\theta}{2}\right) \right\}, \\
F_l(r,\theta) \equiv \sqrt{r} \left\{
\sin\left(\frac{\theta}{2}\right),
\cos\left(\frac{\theta}{2}\right),
\sin\left(\frac{\theta}{2}\right)\sin\left(\theta\right),
\cos\left(\frac{\theta}{2}\right)\sin\left(\theta\right)\right\}.
\end{eqnarray}
\label{eqn:asymp}
\end{subequations}

Here $(r,\theta)$ are polar coordinates in the local coordinate system with the origin at the crack tip. The functions described here to recover the singular fields around the crack tip were originally proposed for isotropic plates~\cite{Dolbow2000}. As we are interested in the global behavior of the cracked FGM plate, we propose to use the same enrichment functions. The role of these enrichment functions is to aid in representing the discontinuous surface independent of the mesh.

In \Erefs{eqn:platexfem1} and (\ref{eqn:platexfem2}), $\mathcal{N}^{\rm{fem}}$ is the set of all nodes in the mesh; $\mathcal{N}^{\rm{c}}$ is the set of nodes whose shape function support is cut by the crack interior (squared nodes in \fref{fig:xfem_crack}) and  $\mathcal{N}^{\rm{f}}$ is the set of nodes whose shape function support is cut by the crack tip (circled nodes in \fref{fig:xfem_crack}). For any node in $\mathcal{N}^{\rm{f}}$, the support of the nodal shape function is fully cut into two disjoint pieces by the crack. 
If for a certain node $n_i$, one of the two pieces is very small compared to the other, then the generalized Heaviside function used for the enrichment is almost a constant over the support, leading to an ill-conditioned matrix~\cite{Moes1999}. Therefore, in this case, the node $n_i$ is removed from the set $\mathcal{N}^{\rm{c}}$. The area-criterion for  the nodal inclusion in $\mathcal{N}^{\rm{c}}$ is as follows: let the area above the crack is $A^{ab}_{\omega}$ and the area below the crack is $A^{be}_\omega$ and $A_\omega = A^{ab}_\omega + A^{be}_\omega$. If either of the two ratios, $A^{ab}_\omega / A_\omega$ or $A^{be}_\omega / A_\omega$ is below a prescribed tolerance, the node is removed from the set $\mathcal{N}^{\rm{c}}$. A tolerance of 10$^{-4}$ is used in the computations.

\subsection{Discretized equations for enriched Q4 plate element}

Now, applying the displacement field approximated
by~\Erefs{eqn:platexfem1} - (\ref{eqn:platexfem2}) in \Eref{eqn:potential} and
\Eref{eqn:kinetic}, one gets the modified Lagrange's equations of
motion in discretized form as:

\begin{equation}
\left( \tilde{\KK} - \omega^2 \tilde{\bm} \right) \boldsymbol{\delta} = \mathbf{0},
\label{eqn:vibeqn1}
\end{equation}

where the element stiffness matrix is given by:

\begin{equation}
\renewcommand\arraystretch{2}
\tilde{\KK}^e = \left[\begin{array}{*{20}c}
               \tilde{\KK}_{uu}^e  & \tilde{\KK}_{ua}^e \\
                 \tilde{\KK}_{au}^e  &  \tilde{\KK}_{aa}^e
                 \end{array} \right] = \int_{\Omega^e}\left[ \begin{array}{*{20}c}
                 \BB_{std}^T \mathbf{D} \BB_{std} &  \BB_{std}^T \mathbf{D}
                 \BB_{enr} \\  \BB_{enr}^T \mathbf{D}
                 \BB_{std} &  \BB_{enr}^T \mathbf{D}
                 \BB_{enr} \end{array} \right]~d\Omega_e,
\end{equation}

where $\BB_{std}$ and $\BB_{enr}$ are the standard and
enriched part of the strain-displacement matrix, respectively and $\mathbf{D}$ is the
material matrix. The element mass matrix is given by:

\begin{equation}
\renewcommand\arraystretch{2}
\tilde{\bm}^e = \left[\begin{array}{*{20}c}
               \tilde{\bm}_{uu}^e  & \tilde{\bm}_{ua}^e \\
                 \tilde{\bm}_{au}^e  &  \tilde{\bm}_{aa}^e
                 \end{array} \right] = \int_{\Omega^e}\left[ \begin{array}{*{20}c}
                 \bn_{std}^T ~\rho \bn_{std} &  \bn_{std}^T ~\rho
                 \bn_{enr} \\  \bn_{enr}^T ~\rho
                 \bn_{std} & \bn_{enr}^T ~\rho
                 \bn_{enr}\end{array} \right] ~d\Omega_e.
\end{equation}

where in deriving the above element mass matrix, the plate displacements and the section rotations given by~\Erefs{eqn:platexfem1} - (\ref{eqn:platexfem2}) are used. In solving for the eigenvalues, the QR algorithm, based on the QR decomposition is used~\cite{Watkins2002}. 

\section{Numerical results}
\label{numerics} 
In this section, we present the natural frequencies of a cracked functionally graded material plates using the extended Q4 formulation. We consider both square and rectangular plates with simply supported, cantilevered and clamped boundary conditions. In all cases, we present the non dimensionalized free flexural frequencies as, unless specified otherwise:

\begin{equation}
\Omega = \omega \left( \frac{b^2}{h} \right) \sqrt{ \frac{\rho_c}{E_c}}
\label{eqn:nondimfreq}
\end{equation}

where $E_c, \nu_c$ are the Young's moulus and Poisson's ratio of the ceramic material, 
and $\rho_c$ is the mass density. In order to be consistent with the existing literature, properties of the ceramic are used for normalization. The effect of plate thickness $a/h$, aspect ratio $b/a$, crack length $d/a$, crack orientation $\theta$, location of the crack, multiple cracks and boundary condition on the natural frequencies are studied. Based on progressive mesh refinement, a $34\times34$ structured mesh is found to be adequate to model the full plate for the present analysis. The material properties used for the FGM components are listed in Table~\ref{table:matprop}. 

Before proceeding with the detailed study on the effect of different parameters on the natural frequency, the formulation developed herein is validated against available results pertaining to the linear frequencies of cracked isotropic and functionally graded material plates with different boundary conditions.  The computed frequencies for cracked isotropic simply supported rectangular plate is given in  Table~\ref{table:RectComparison}. Tables~\ref{table:SSconvresults} and \ref{table:Canticonvresults} gives a comparison of computed frequencies for simply supported square plate with a side crack and cantilevered plate with a side crack, respectively. It can be seen that the numerical results from the present formulation are found to be in  good agreement with the existing solutions.

The FGM plate considered here consists of silicon nitride (SI$_3$N$_4$) and stainless steel (SUS304). The material is considered to be temperature dependent and the temperature coefficients corresponding to SI$_3$N$_4$/SUS304 are listed in Table \ref{table:tempdepprop} ~\cite{Sundararajan2005,Reddy1998}. The mass density $(\rho)$ and thermal conductivity $(K)$ are: $\rho_c$=2370 kg/m$^3$, $K_c$=9.19 W/mK for SI$_3$N$_4$ and $\rho_m$ = 8166 kg/m$^3$, $K_m$ = 12.04 W/mK for SUS304. Poisson's ratio $\nu$ is assumed to be constant and taken as 0.28 for the current study~\cite{Sundararajan2005,Prakash2007}.  Here, the modified shear correction factor obtained based on energy equivalence principle as outlined in~\cite{Singh2011} is used.
The boundary conditions for simply supported and clamped cases are (see \fref{fig:SS2Plate}):

\noindent \emph{Simply supported boundary condition}: \\
\begin{eqnarray}
u_o = w_o = \theta_y = 0 \hspace{1cm} ~\textup{on} ~ x=0,a \nonumber \\
v_o = w_o = \theta_x = 0 \hspace{1cm} ~\textup{on} ~ y=0,b
\end{eqnarray}

\noindent \emph{Clamped boundary condition}: \\
\begin{eqnarray}
u_o = w_o = \theta_y = v_o = \theta_x = 0 \hspace{1cm} ~\textup{on} ~ x=0,a \nonumber \\
u_o = w_o = \theta_y = v_o = \theta_x = 0 \hspace{1cm} ~\textup{on} ~ y=0,b
\end{eqnarray}

\begin{table}[htpb]
\renewcommand\arraystretch{1.5}
\caption{Material properties of the FGM components. $^\dag$Ref~\cite{Huang2011}, $^\ast$\cite{Reddy1998,Sundararajan2005}}
\centering
\begin{tabular}{lccc}
\hline
Material & \multicolumn{3}{c}{Properties} \\
\cline{2-4}
& E(GPa) & $\nu$ & $\rho$ (Kg/m$^3$) \\
\hline 
Aluminum (Al)$^\dag$ & 70.0 & 0.30 & 2702 \\
Alumina (Al$_2$O$_3$)$^\dag$ & 380.0 & 0.30 & 3800 \\
Zirconia (ZrO$_2$)$^\dag$ & 200.0 & 0.30 & 5700 \\
Steel (SUS304)$^\ast$ & 201.04 & 0.28 & 8166 \\
Silicon Nitride (Si$_3$N$_4$)$^\ast$ & 348.43 & 0.28 & 2370 \\
\hline 
\end{tabular}
\label{table:matprop}
\end{table}

\begin{table}
\renewcommand\arraystretch{1.5}
\caption{Temperature dependent coefficient for material SI$_3$N$_4$/SUS304, Ref~\cite{Reddy1998,Sundararajan2005}.}
\centering
\begin{tabular}{lcccccc}
\hline
Material & Property & $P_o$ & $P_{-1}$ & $P_1$ & $P_2$ & $P_3$  \\
\hline
\multirow{2}{*}{SI$_3$N$_4$} & $E$(Pa) & 348.43e$^9$ &0.0& -3.070e$^{-4}$ & 2.160e$^{-7}$ & -8.946$e^{-11}$  \\
& $\alpha$ (1/K) & 5.8723e$^{-6}$ & 0.0 & 9.095e$^{-4}$ & 0.0 & 0.0 \\
\cline{2-7}
\multirow{2}{*}{SUS304} & $E$(Pa) & 201.04e$^9$ &0.0& 3.079e$^{-4}$ & -6.534e$^{-7}$ & 0.0  \\
& $\alpha$ (1/K) & 12.330e$^{-6}$ & 0.0 & 8.086e$^{-4}$ & 0.0 & 0.0 \\
\hline
\end{tabular}
\label{table:tempdepprop}
\end{table}

\begin{table}[htpb]
\renewcommand\arraystretch{1.5}
\caption{Comparison of frequency parameters $\omega(b^2/h)\sqrt{\rho_c/E_c}$ for a simply supported homogeneous rectangular thin plate with a horizontal crack $(a/b = 2, b/h = 100, c_y/b = 0.5, d/a = 0.5, \theta = 0)$.}
\centering
\begin{tabular}{ccccc}
\hline 
mode & Ref~\cite{Stahl1972} & Ref~\cite{Huang2009} & Ref~\cite{Huang2011} & Present \\
\hline
1 & 3.050 & 3.053 & 3.047 & 3.055 \\
2 & 5.507 & 5.506 & 5.503 & 5.508 \\
3 & 5.570 & 5.570 & 5.557 & 5.665 \\
4 & 9.336 & 9.336 & 9.329 & 9.382 \\
5 & 12.760 & 12.780 & 12.760 & 12.861\\
\hline 
\end{tabular}
\label{table:RectComparison}
\end{table}

\begin{table}[htpb]
\renewcommand\arraystretch{1.5}
\caption{Non-dimensionalized natural frequency for a simply supported
square Al/Al$_2$O$_3$ plate with a side crack $(a/b=1, a/h=10)$. Crack length $d/a=0.5$.} \centering
\begin{tabular}{crrrrrr}
\hline
gradient & \multicolumn{2}{c}{Mode 1} & \multicolumn{2}{c}{Mode 2} & \multicolumn{2}{c}{Mode 3}\\
\cline{2-7}
index, $n$ & Ref~\cite{Huang2011}  & Present & Ref~\cite{Huang2011}  & Present & Ref~\cite{Huang2011}  & Present \\
\hline
0 & 5.379 & 5.387 & 11.450 & 11.419 & 13.320 & 13.359 \\
0.2 & 5.001 & 5.028 & 10.680 & 10.659 & 12.410 & 12.437 \\
1 & 4.122 & 4.122 & 8.856 & 8.526 & 10.250 & 10.285 \\
5 & 3.511 & 3.626 & 7.379 & 7.415 & 8.621 & 8.566 \\
10 & 3.388 & 3.409 & 7.062 & 7.059 & 8.289 & 8.221 \\
\hline
\end{tabular}
\label{table:SSconvresults}
\end{table}

\begin{table}[htpb]
\renewcommand\arraystretch{1.5}
\caption{Fundamental frequency $\omega b^2/h \sqrt{\rho_c/E_c}$ for cantilevered square Al/ZrO$_2$ FGM plates with horizontal size crack $(b/h=10, c_y/b = 0.5, d/a = 0.5)$.} \centering
\begin{tabular}{cccccccc}
\hline
$a/b$& Mode & & \multicolumn{5}{c}{gradient index, $n$} \\
\cline{4-8}
& & & 0 & 0.2 & 1 & 5 & 10 \\
\hline
\multirow{6}{*}{1} & \multirow{2}{*}{1} & Ref~\cite{Huang2011} & 1.0380 & 1.0080 & 0.9549 & 0.9743 & 0.9722 \\
& & Present & 1.0380 & 1.0075 & 0.9546 & 0.9748 & 0.9722 \\
& \multirow{2}{*}{2} & Ref~\cite{Huang2011} & 1.7330 & 1.6840 & 1.5970 & 1.6210 & 1.6170 \\
& & Present & 1.7329 & 1.6834 & 1.5964 & 1.6242 & 1.6194 \\
& \multirow{2}{*}{3} & Ref~\cite{Huang2011} & 4.8100 & 4.6790 & 4.4410 & 4.4760 & 4.4620 \\
& & Present & 4.8231 & 4.6890 & 4.4410 & 4.4955 & 4.4845 \\
\hline
\end{tabular}
\label{table:Canticonvresults}
\end{table}

\input{SS2_ah10_results}

\input{SS2_ah10_multicrk}

\subsection*{Effect of aspect ratio, thickness and boundary conditions}
The influence of plate aspect ratio $b/a$, plate thickness $a/h$ and boundary condition on a cracked FGM plate with a horizontal center crack is shown in Table~\ref{table:FreqAhAb}. Two types of boundary conditions are studied: Simply supported (SS) and Clamped condition (CC). For a given crack length and for a given crack location, decreasing the plate thickness and increasing the plate aspect ratio, increases the frequency. The increase in stiffness is the cause for increase in frequency when the boundary condition is changed from SS to CC for a fixed aspect ratio and plate thickness. 

\begin{table}[htpb]
\renewcommand\arraystretch{1.2}
\caption{Effect of plate aspect ratio $b/a$, plate thickness $a/h$ and boundary condition on fundamental frequency $\omega b^2/h \sqrt{\rho_c/E_c}$ for Si$_3$N$_4$/SUS304 FGM plate with horizontal center crack $(c_y/b = 0.5, d/a = 0.5)$. $^\dag$Simply Supported, $^\ast$Clamped Support.}
\centering
\begin{tabular}{ccccccc}
\hline
$b/a$ & $a/h$ & \multicolumn{2}{c}{Mode 1} & & \multicolumn{2}{c}{Mode 2} \\
\cline{3-4} \cline{6-7}
& & SS$^\dag$ & CC$^\ast$ & & SS$^\dag$ & CC$^\ast$  \\
\hline 
\multirow{3}{*}{0.5} & 10 & 1.1205 & 2.2202 && 2.1586 & 2.8588 \\
& 20 & 1.1974 & 2.5043 && 2.5482 & 3.4621  \\
& 100 & 1.2625 & 2.6748 && 2.6593 & 4.0903 \\
\cline{2-7}
\multirow{3}{*}{1} & 10 & 2.4051 & 4.1624 && 5.2792 & 6.6286 \\
& 20 & 2.4831 & 4.4592 && 5.8338 & 7.6464 \\
& 100 & 2.5473 & 4.6311 && 6.2765 & 8.5774 \\
\cline{2-7}
\multirow{3}{*}{2} & 10 & 6.6864 & 12.7513 && 10.5295 & 15.9115 \\
& 20 & 6.8101 & 13.4551 && 10.8647 & 17.0301 \\
& 100 & 6.8847 & 13.7540 &&11.0392 & 17.6030 \\
\hline
\end{tabular}
\label{table:FreqAhAb}
\end{table}

\subsection{Plate with side crack}
\begin{figure}[htpb]
\centering
\subfigure[Effect of aspect ratio, $a/h$]{\includegraphics[scale=0.6]{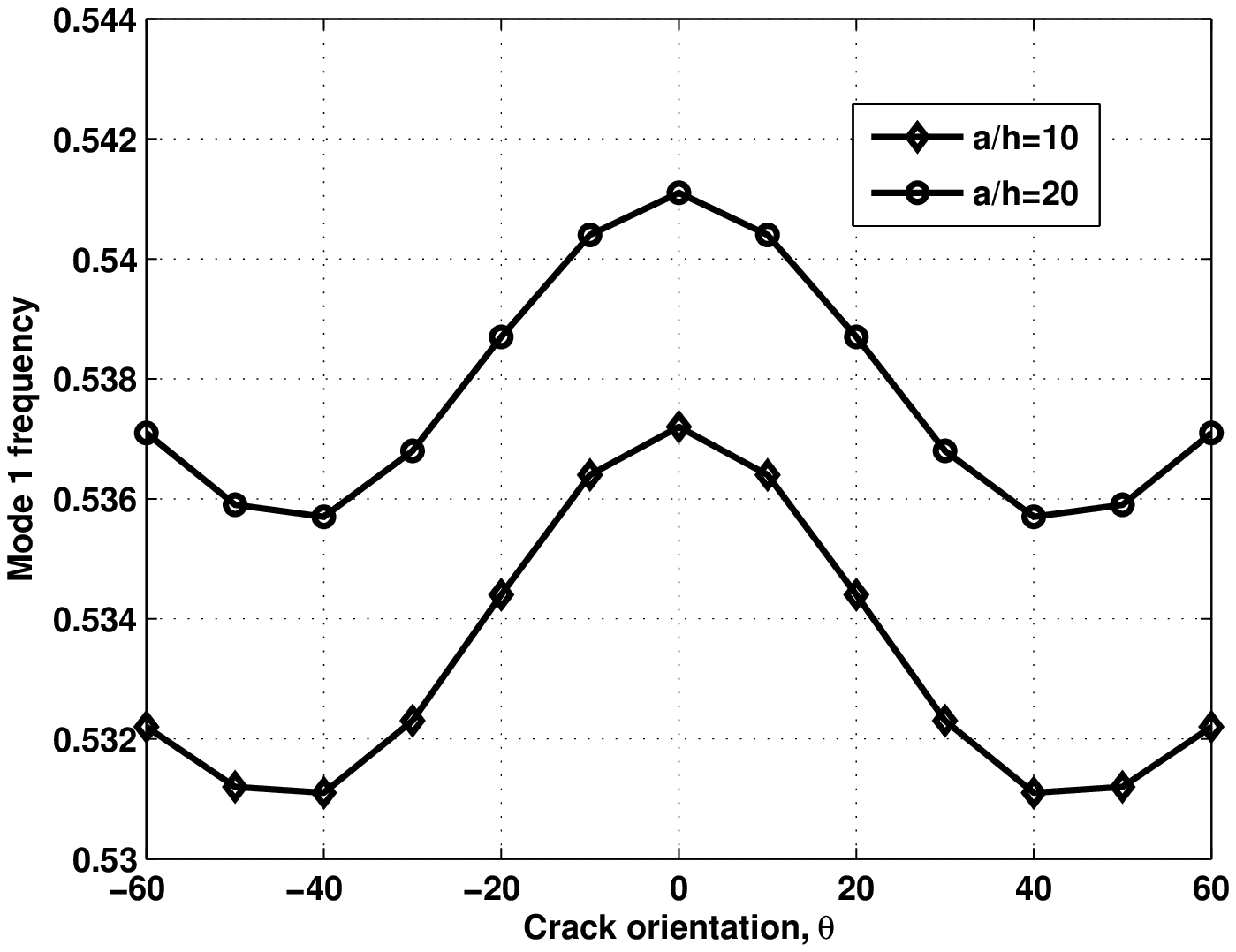}}
\subfigure[Effect of gradient index, $n$]{\includegraphics[scale=0.6]{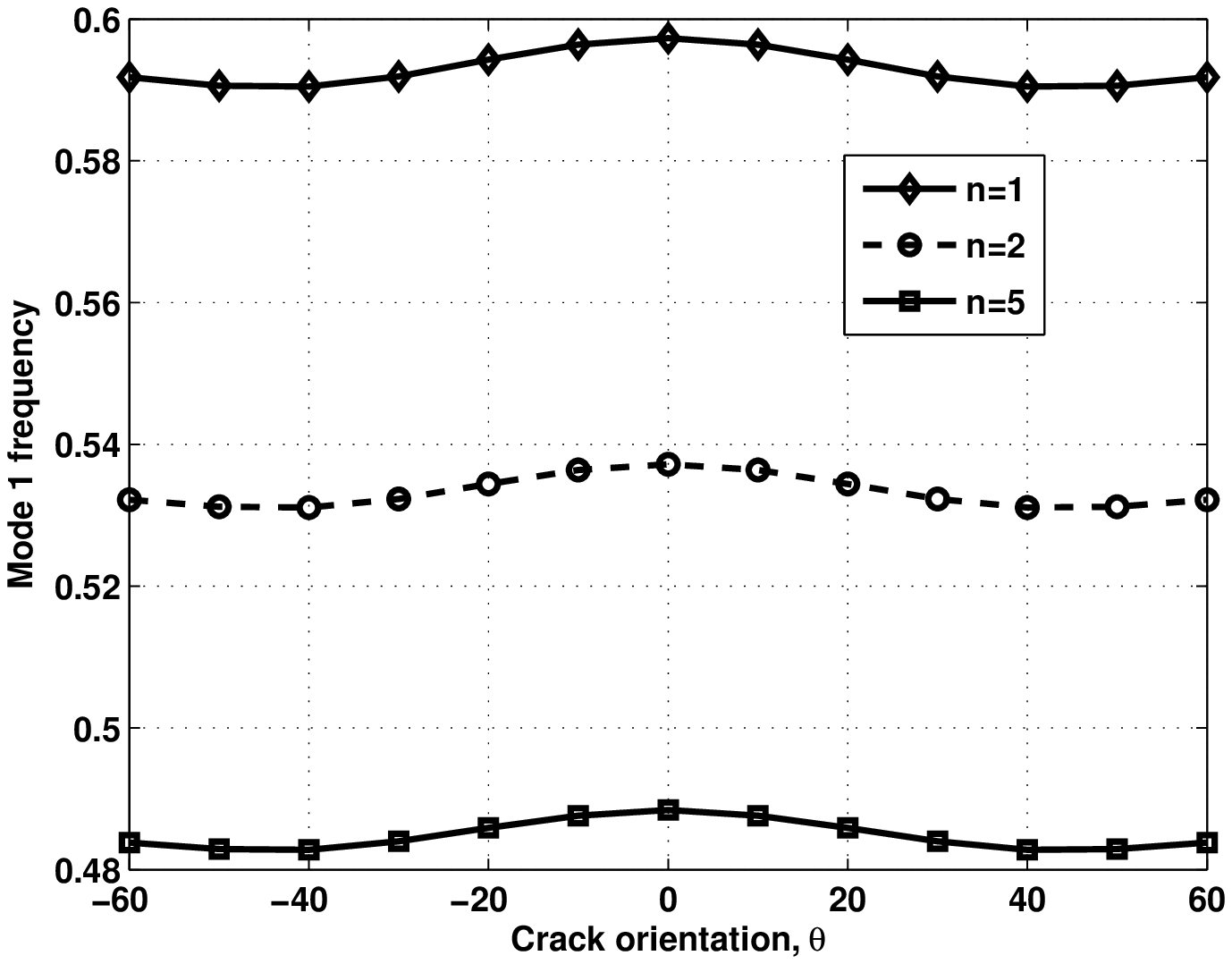}}
\caption{Variation of fundamental frequency for a cantilevered square Si$_3$N$_4$/SUS304 FGM plate with a side crack as a function of crack orientation. $d/a=$0.5, $c_y/a$=0.5.}
\label{fig:Canti}
\end{figure}

\begin{figure}[htpb]
\centering
\subfigure[Without crack, $\omega_1=$0.4885]{\includegraphics[scale=0.6]{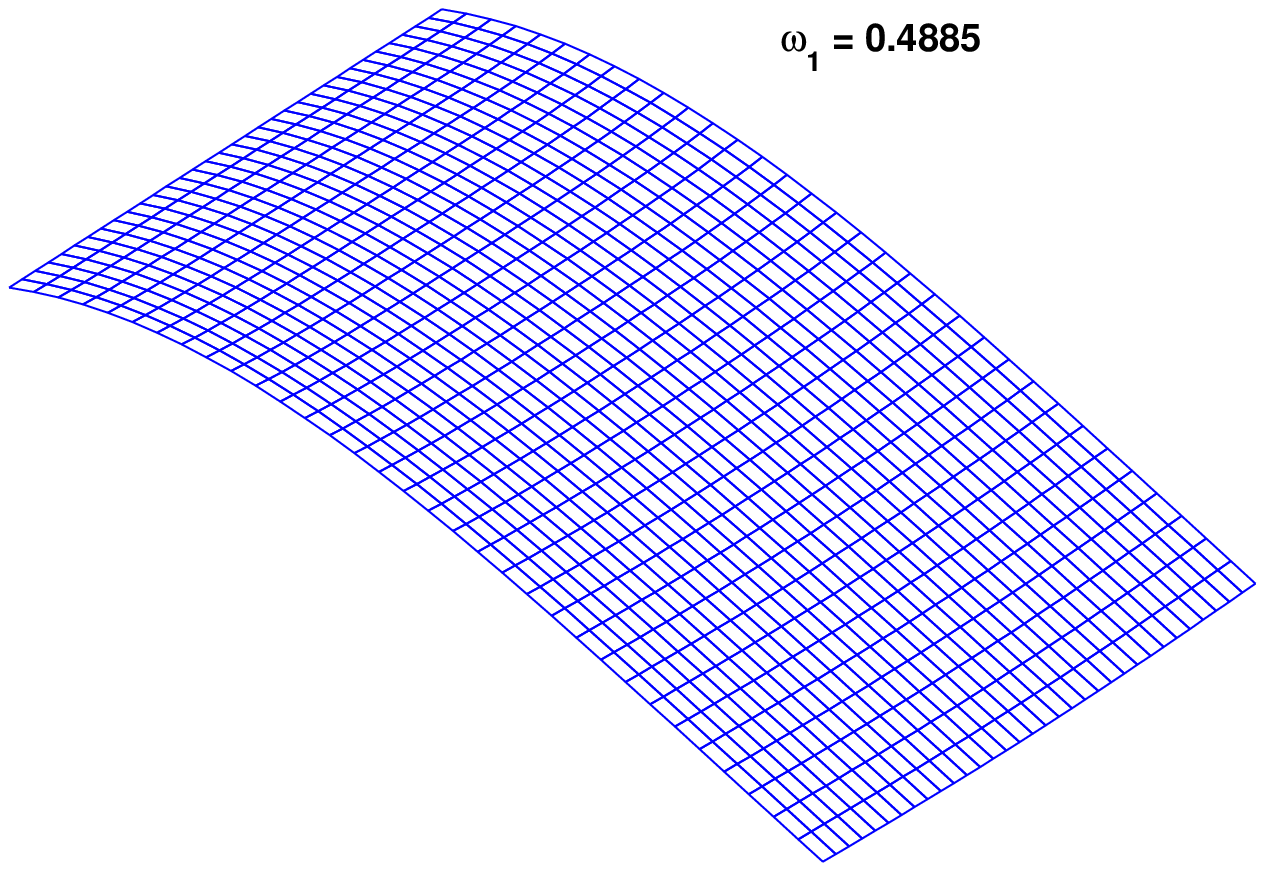}}
\subfigure[With crack, $\omega_1=$0.4883]{\includegraphics[scale=0.6]{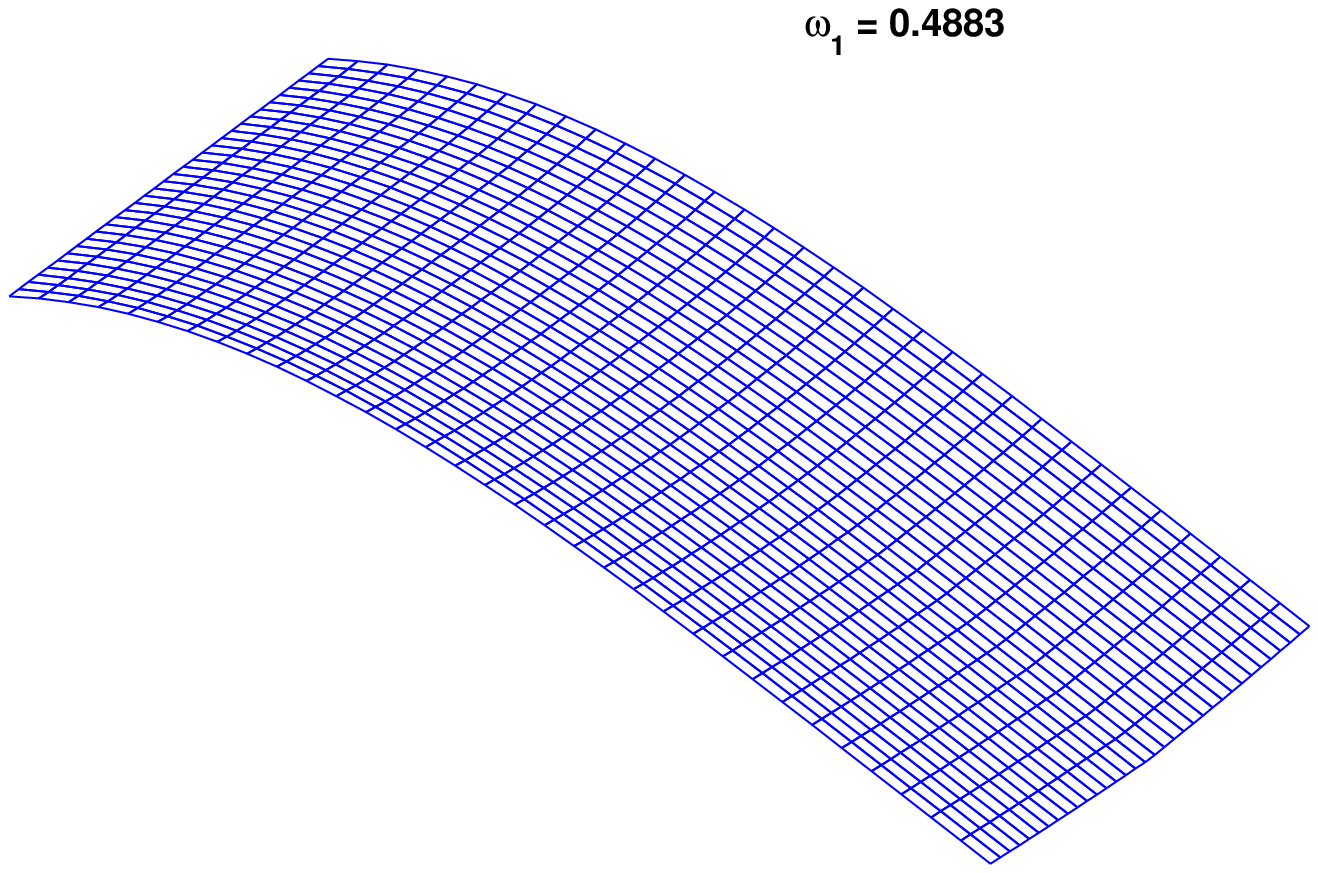}}
\caption{First Mode shape for a cantilevered plate with a side crack with $\theta = 0, n = 5, c_y/b = 0.5, d/a = 0.5, a/h=10, b/a=1$.}
\label{fig:cantimode1shapes}
\end{figure}

\begin{figure}[htpb]
\centering
\subfigure[Without crack, $\omega_2=$1.1608]{\includegraphics[scale=0.6]{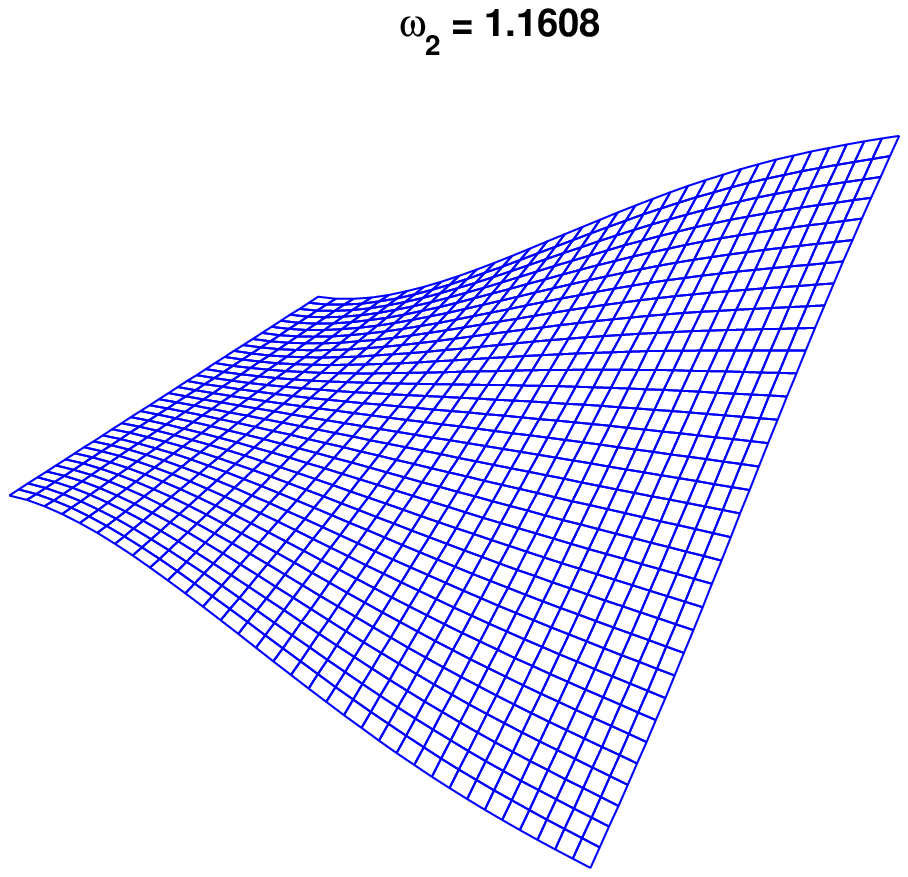}}
\subfigure[With crack, $\omega_2=$0.8223]{\includegraphics[scale=0.6]{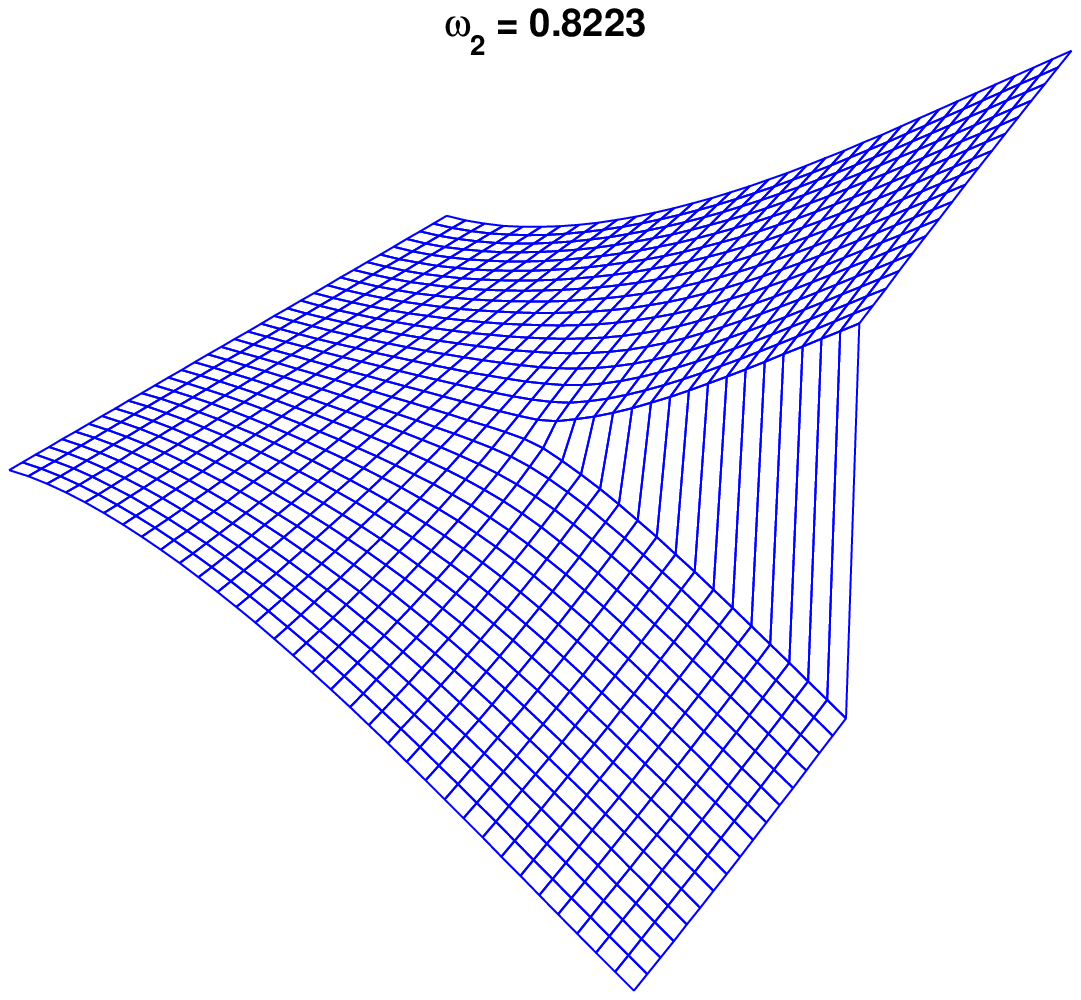}}
\caption{Second Mode shape for a cantilevered plate with a side crack with $\theta = 0, n = 5, c_y/b = 0.5, d/a = 0.5, a/h=10, b/a=1$.}
\label{fig:cantimode2shapes}
\end{figure}

\begin{figure}[htpb]
\centering
\subfigure[$\theta=40^o$]{\includegraphics[scale=0.6]{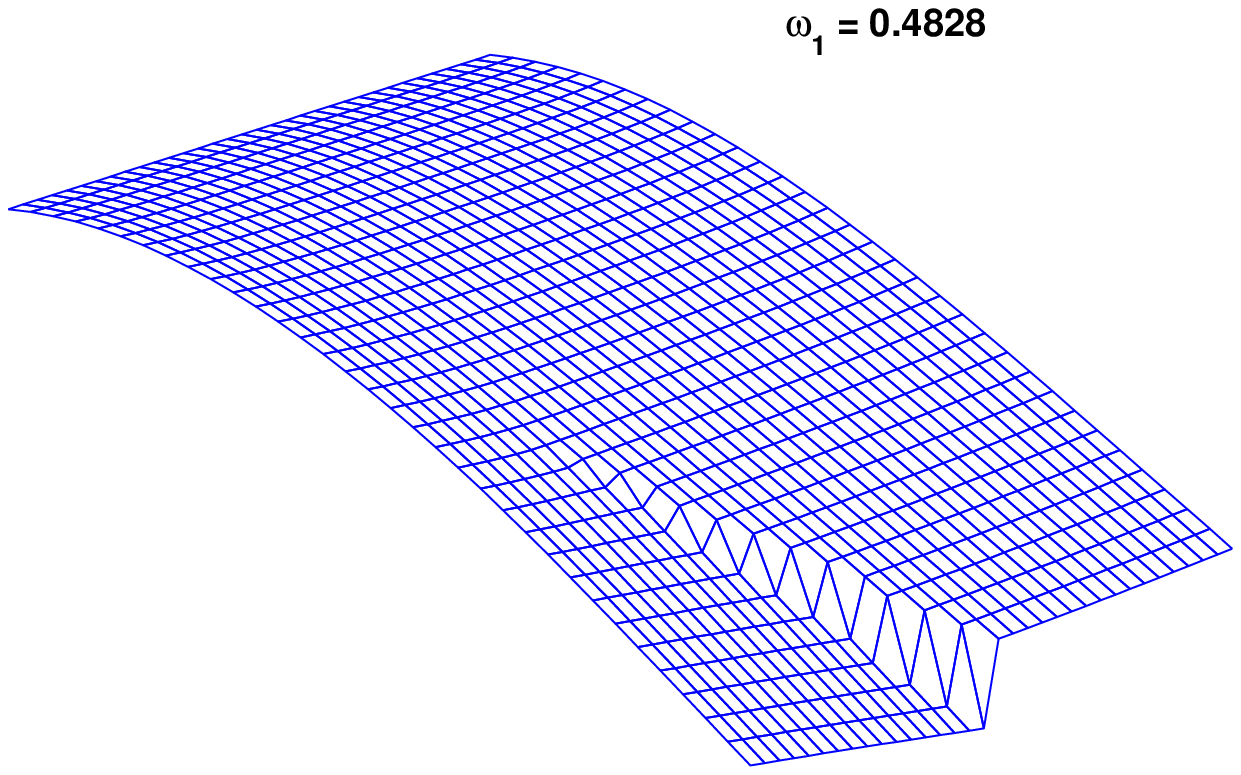}}
\subfigure[$\theta=-40^o$]{\includegraphics[scale=0.6]{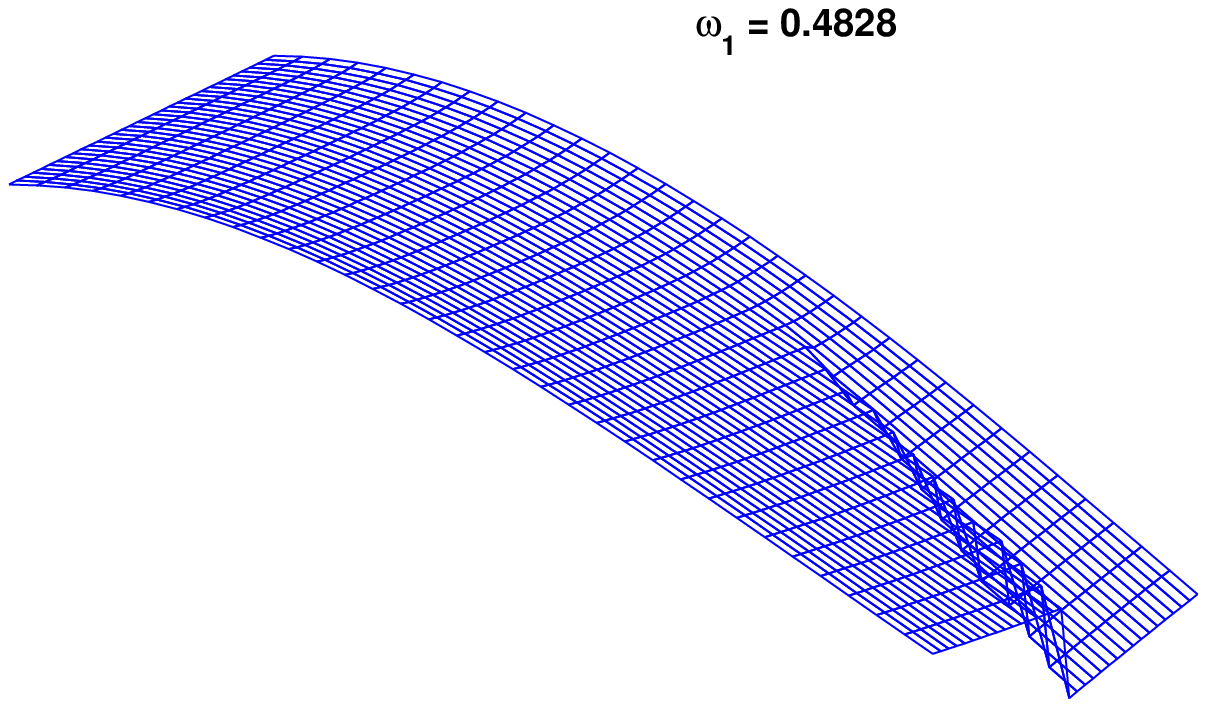}}
\caption{First Mode shape $(\omega_1=0.4828)$ for a cantilevered plate with a side crack with $n = 5, c_y/b = 0.5, d/a = 0.5, a/h=10, b/a=1$.}
\label{fig:cantimodeangle40}
\end{figure}

\begin{figure}[htpb]
\centering
\subfigure[$\theta=60^o$]{\includegraphics[scale=0.6]{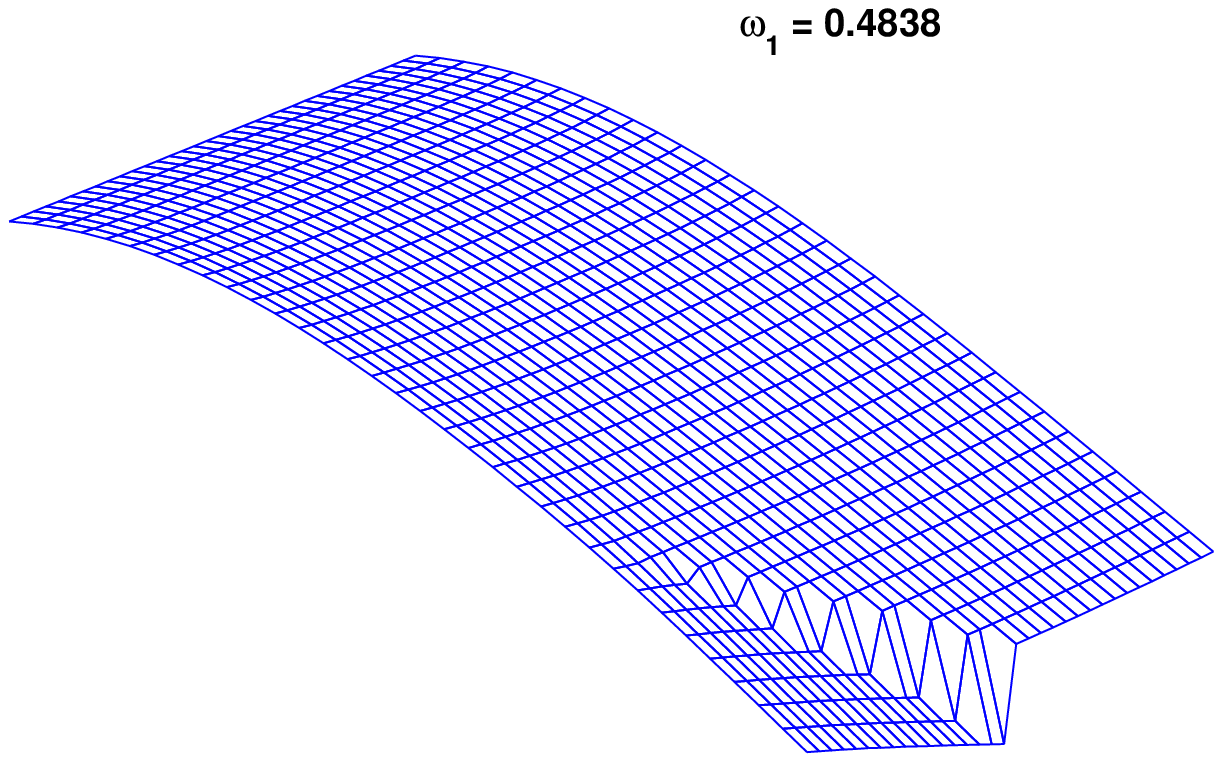}}
\subfigure[$\theta=-60^o$]{\includegraphics[scale=0.6]{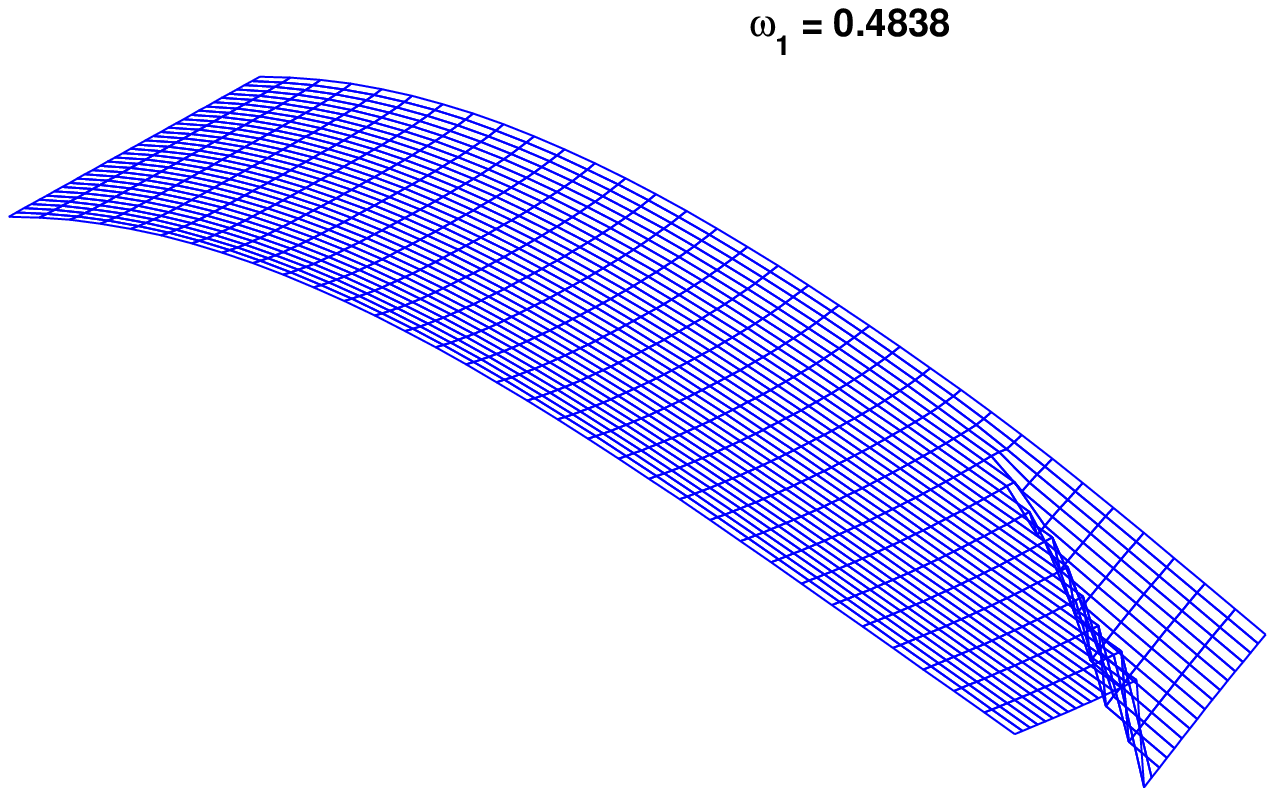}}
\caption{First Mode shape $(\omega_1=0.4838)$ for a cantilevered plate with a side crack with $n = 5, c_y/b = 0.5, d/a = 0.5, a/h=10, b/a=1$.}
\label{fig:cantimodeangle60}
\end{figure}

Consider a plate of uniform thickness, $h$, with length and width as $a$ and $b$, respectively. \fref{fig:cantiplate} shows a cantilevered plate with a side crack of length $d$ located at a distance of $c_y$ from the x-axis and at an angle $\theta$ with respect to the $x-$ axis. The influence of plate thickness, crack orientation and gradient index on the fundamental frequency is shown in Table~\ref{table:Cantiaspect} and in \fref{fig:Canti}. With increase in gradient index $n$, the frequency decreases for all crack orientations and for different plate thickness. With increase in crack orientation from $\theta = -60^o$ to $\theta = 60^o$, the frequency initially decreases until $\theta = -40^o$ and then reaches the maximum when the crack is horizontal. And with further increase in the crack orientation, the frequency decreases and the plate's response is symmetric. This is because, when the crack is horizontal $(\theta=0)$, the crack is aligned to the first mode shape and the response of the plate is similar to the cantilevered plate without a crack. The frequency of the plate without a crack is greater than a plate with a crack. \frefs{fig:cantimode1shapes} and (\ref{fig:cantimode2shapes}) shows first two mode shapes for a cantilevered plate with and without a horizontal crack. As explained earlier the frequency and the first mode shape for a plate with and without a crack are very similar. For any other crack orientation, the mode shape would be influenced by the presence of the crack. \frefs{fig:cantimodeangle40} and (\ref{fig:cantimodeangle60}) shows the first mode shape of a cantilevered plate with a side crack with orientations $\theta=\pm 40^o$ and $\theta = \pm 60^o$, respectively.

\begin{figure}[htpb]
\centering
\input{cantilever.pstex_t}
\caption{Cantilevered plate with a side crack: geometry}
\label{fig:cantiplate}
\end{figure}

\begin{table}[htpb]
\renewcommand\arraystretch{1.2}
\caption{Fundamental frequency $\omega b^2/h \sqrt{\rho_c/E_c}$ for cantilevered square plate Si$_3$N$_4$/SUS304 FGM plate with a side crack $(c_y/b = 0.5, d/a = 0.5)$ as a function of crack angle and gradient index. $^\dag$denotes change in trend.}
\centering
\begin{tabular}{ccccccc}
\hline
$a/h$& crack & \multicolumn{5}{c}{gradient index, $n$}\\
\cline{3-7}
& angle, $\theta$ & 0 & 1 & 2 & 5 & 10 \\
\hline
\multirow{13}{*}{10} &-60	&0.9859	&0.5918	&0.5322	&0.4838	&0.4610\\
&-50	&0.9840	&0.5906	&0.5312	&0.4829	&0.4601\\
&{\bf -40}$^\dag$ & {\bf 0.9838} & {\bf 0.5905 } & {\bf 0.5311} & {\bf 0.4828} & {\bf 0.4600}\\
&-30	&0.9862	&0.5919	&0.5323	&0.4840	&0.4611\\
&-20	&0.9900	&0.5943	&0.5344	&0.4859	&0.4630\\
&-10	&0.9936	&0.5964	&0.5364	&0.4876	&0.4646\\
&{\bf 0}$^\dag$ & {\bf 0.9951}	& {\bf 0.5973} & {\bf 0.5372} &{\bf 0.4884} &{\bf 0.4653}\\
&10	&0.9936	&0.5964	&0.5364	&0.4876	&0.4646\\
&20	&0.9900	&0.5943	&0.5344	&0.4859	&0.4630\\
&30	&0.9862	&0.5919	&0.5323	&0.4840	&0.4611\\
&{\bf 40}$^\dag$ & {\bf 0.9838} & {\bf 0.5905} & {\bf 0.5311} & {\bf 0.4828} & {\bf 0.4600}\\
&50	&0.9840	&0.5906	&0.5312	&0.4829	&0.4601\\
&60	&0.9859	&0.5918	&0.5322	&0.4838	&0.4610\\
\cline{2-7}
\multirow{13}{*}{20} &-60 &0.9949 &0.5972 &0.5371 &0.4883 &0.4653\\
&-50 &0.9927 &0.5959 &0.5359 &0.4872 &0.4643\\
&{\bf -40}$^\dag$ &{\bf 0.9924} & {\bf 0.5957} & {\bf 0.5357} & {\bf 0.4871} & {\bf 0.4641}\\
&-30 &0.9944 &0.5969 &0.5368 &0.4881 &0.4651\\
&-20 &0.9979 &0.5989 &0.5387 &0.4898 &0.4667\\
&-10 &1.0011 &0.6009 &0.5404 &0.4913 &0.4682\\
&{\bf 0}$^\dag$ &{\bf 1.0024} &{\bf 0.6016} &{\bf 0.5411} &{\bf 0.4920} &{\bf 0.4688}\\
&10	&1.0011	&0.6009	&0.5404	&0.4913	&0.4682\\
&20	&0.9979 &0.5989	&0.5387	&0.4898	&0.4667\\
&30	&0.9944	&0.5969	&0.5368	&0.4881	&0.4651\\
&{\bf 40}$^\dag$ &{\bf 0.9924} & {\bf 0.5957} & {\bf 0.5357} & {\bf 0.4871} & {\bf 0.4641}\\
&50	&0.9927	&0.5959	&0.5359	&0.4872	&0.4643\\
&60 &0.9949	&0.5972	&0.5371	&0.4883	&0.4653\\
\hline
\end{tabular}
\label{table:Cantiaspect}
\end{table}

\section{Conclusion}
\label{concl}
Natural frequencies of cracked functionally graded material plate is studied using the extended finite element method. The formulation is based on first order shear deformation theory for plates and a four-noded field consistent enriched element is used. The material is assumed to be temperature dependent and graded in the thickness direction. Numerical experiments have been conducted to bring out the effect of gradient index, crack length, crack orientation, crack location, boundary condition, plate aspect ratio and plate thickness on the natural frequency of the FGM plate. Also the influence of multiple cracks and their relative orientation on the natural frequency is studied. From the detailed numerical study, the following can be concluded:

\begin{itemize}
\item Increasing the crack length decreases the natural frequency. This is due to the reduction in stiffness of the material structure. The frequency is lowest when the crack is located at the center of the plate. 
\item Increasing gradient index $n$ decreases the natural frequency. This is due to the increase in metallic volume fraction.
\item Decreasing the plate thickness $a/h$ and increasing the aspect ratio $b/a$ increases the frequency. 
\item For a cantilevered plate with a side crack, horizontal crack has the maximum frequency and the trend changes at $\theta = \pm 40^o$.
\item Crack orientation $\theta = 45^o$ has been observed to be a critical angle. At this crack orientation, the frequency changes its trend for a square plate.
\item Increasing the number of cracks, decreases the overall stiffness of the plate and thus decreasing the frequency. The frequency is lowest when both the cracks are oriented at $\theta = 50^o$.
\end{itemize}

{\bf Acknowledgements} \\

S Natarajan acknowledges the financial support of (1) the Overseas Research Students Awards Scheme; (2) the Faculty of Engineering, University of Glasgow, for period Jan. 2008 - Sept. 2009 and of (3) the School of Engineering (Cardiff University) for the period Sept. 2009 onwards.

S Bordas would like to thank the support of the Royal Academy of Engineering and of the Leverhulme Trust for his Senior Research Fellowship entitled ``Towards the next generation surgical simulators" as well as the support of EPSRC under grants EP/G042705/1 Increased Reliability for Industrially Relevant Automatic Crack Growth Simulation with the eXtended Finite Element Method.

\bibliographystyle{elsarticle-num}
\bibliography{myVibCrk}

\end{document}

%% file: plate.pstex_t
\begin{picture}(0,0)%
\includegraphics{plate.pstex}%
\end{picture}%
\setlength{\unitlength}{2279sp}%
\begingroup\makeatletter\ifx\SetFigFont\undefined%
\gdef\SetFigFont#1#2#3#4#5{%
  \reset@font\fontsize{#1}{#2pt}%
  \fontfamily{#3}\fontseries{#4}\fontshape{#5}%
  \selectfont}%
\fi\endgroup%
\begin{picture}(7002,4165)(2689,-5264)
\put(4501,-5191){\makebox(0,0)[lb]{\smash{{\SetFigFont{8}{9.6}{\familydefault}{\mddefault}{\updefault}$a$}}}}
\put(8371,-3481){\makebox(0,0)[lb]{\smash{{\SetFigFont{8}{9.6}{\familydefault}{\mddefault}{\updefault}$b$}}}}
\put(9676,-1861){\makebox(0,0)[lb]{\smash{{\SetFigFont{8}{9.6}{\familydefault}{\mddefault}{\updefault}$h$}}}}
\end{picture}%

%% file: XFEMElemcate.pstex_t
\begin{picture}(0,0)%
\includegraphics{XFEMElemcate.pstex}%
\end{picture}%
\setlength{\unitlength}{2486sp}%
\begingroup\makeatletter\ifx\SetFigFont\undefined%
\gdef\SetFigFont#1#2#3#4#5{%
  \reset@font\fontsize{#1}{#2pt}%
  \fontfamily{#3}\fontseries{#4}\fontshape{#5}%
  \selectfont}%
\fi\endgroup%
\begin{picture}(9653,9099)(1339,-8533)
\put(3781,-6331){\makebox(0,0)[lb]{\smash{{\SetFigFont{9}{10.8}{\familydefault}{\mddefault}{\updefault}$\Xi_{\alpha}$ enriched element, (tip element)}}}}
\put(3766,-6961){\makebox(0,0)[lb]{\smash{{\SetFigFont{9}{10.8}{\familydefault}{\mddefault}{\updefault}$\vartheta$ enriched element (split element)}}}}
\put(3781,-7666){\makebox(0,0)[lb]{\smash{{\SetFigFont{9}{10.8}{\familydefault}{\mddefault}{\updefault}Partially enriched element (blending element)}}}}
\put(3766,-8371){\makebox(0,0)[lb]{\smash{{\SetFigFont{9}{10.8}{\familydefault}{\mddefault}{\updefault}Standard element}}}}
\put(9586,-1231){\makebox(0,0)[lb]{\smash{{\SetFigFont{10}{12.0}{\familydefault}{\mddefault}{\updefault}$J \in \mathcal{N}^{\rm c}$}}}}
\put(9556,-2146){\makebox(0,0)[lb]{\smash{{\SetFigFont{10}{12.0}{\familydefault}{\mddefault}{\updefault}$K \in \mathcal{N}^{\rm f}$}}}}
\put(8731,-6571){\makebox(0,0)[lb]{\smash{{\SetFigFont{9}{10.8}{\familydefault}{\mddefault}{\updefault}Reproducing elements}}}}
\end{picture}%

%% file: SS2_ah10_results.tex
\subsection{Plate with a center crack}

Consider a plate of uniform thickness, $h$ and with length and width as $a$ and $b$, respectively. \fref{fig:SS2Plate} shows a plate with all edges simply supported with a center crack of length $c$. 

\begin{figure}[htpb]
\centering
\input{simplyS.pstex_t}
\caption{Simply supported plate with a center crack.}
\label{fig:SS2Plate}
\end{figure}

\subsection*{Effect of crack length, crack orientation and gradient index}
\begin{figure}[htpb]
\centering
\includegraphics[scale=0.6]{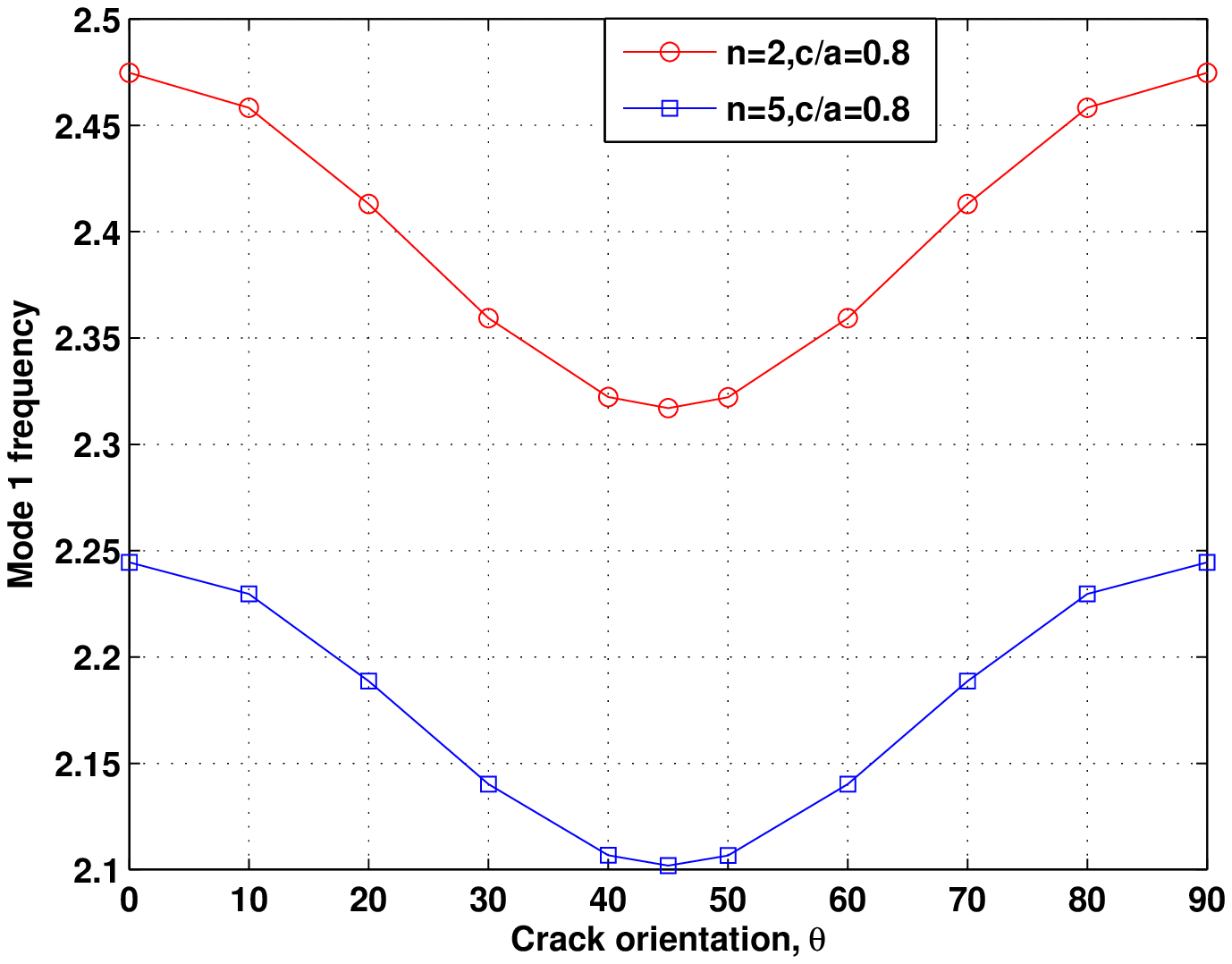}
\caption{Varition of fundamental frequency with orientation of the crack for a simply suppored square FGM plate, $a/h=10$.}
\label{fig:SS2CentCrack}
\end{figure}

The influence of the crack length $d/a$, crack orientation $\theta$ and gradient index $n$ on the fundamental frequency for a simply supported square FGM plate with thickness $a/h=$10 is shown in Tables~\ref{table:SS2ah10k01} and \ref{table:SS2ah10k25}. It is observed that as the crack length increases, the frequency decreases. This is due to the fact that increasing the crack length increases local flexibility and thus decreases the frequency. Also, with increase in gradient index $n$, the frequency decreases. This is because of the stiffness degradation due to increase in metallic volume fraction. It can be seen that the combined effect of increasing the crack length and gradient index is to lower the fundamental frequency. Further, it is observed that the frequency is lowest for a crack orientation $\theta = 45^o$. The frequency values tend to be symmetric with respect to a crack orientation $\theta = 45^o$. This is also shown in \fref{fig:SS2CentCrack} for gradient index $n=5$ and crack length $d/a = 0.8$.

\begin{table}[htpb]
\renewcommand\arraystretch{1.2}
\caption{Fundamental frequency $\omega b^2/h \sqrt{\rho_c/E_c}$ for simply supported Si$_3$N$_4$/SUS304 FGM square plate. $^\dag$denotes change in trend.} \centering
\begin{tabular}{cccccc}
\hline
gradient& Crack & \multicolumn{4}{c}{Crack length, $d/a$.} \\
\cline{3-6}
index, $n$ & orientation, $\theta$ & 0 & 0.4 & 0.6 & 0.8 \\
\hline
\multirow{11}{*}{0}&	0	&	5.5346	&	5.0502	&	4.7526	&	4.5636	\\
&	10	&	5.5346	&	5.0453	&	4.7386	&	4.5337	\\
&	20	&	5.5346	&	5.0379	&	4.7043	&	4.4509	\\
&	30	&	5.5346	&	5.0278	&	4.6640	&	4.3528	\\
&	40	&	5.5346	&	5.0207	&	4.6370	&	4.2849	\\
&	{\bf 45}$^\dag$	&	{\bf 5.5346} & {\bf 5.0173} & {\bf 4.6342} & {\bf 4.2754} \\
&	50	&	5.5346	&	5.0204	&	4.6370	&	4.2849	\\
&	60	&	5.5346	&	5.0278	&	4.6640	&	4.3528	\\
&	70	&	5.5346	&	5.0380	&	4.7043	&	4.4509	\\
&	80	&	5.5346	&	5.0453	&	4.7384	&	4.5337	\\
&	90	&	5.5346	&	5.0503	&	4.7527	&	4.5636	\\
\cline{2-6}
\multirow{11}{*}{1}&	0	&	3.3376	&	3.0452	&	2.8657	&	2.7518	\\
&	10	&	3.3376	&	3.0422	&	2.8571	&	2.7337	\\
&	20	&	3.3376	&	3.0376	&	2.8362	&	2.6833	\\
&	30	&	3.3376	&	3.0315	&	2.8117	&	2.6237	\\
&	40	&	3.3376	&	3.0271	&	2.7953	&	2.5825	\\
&{\bf 45}$^\dag$ & {\bf 3.3376} &{\bf 3.0252} & {\bf 2.7936} &{\bf 2.5767}	\\
&	50	&	3.3376	&	3.0270	&	2.7953	&	2.5824	\\
&	60	&	3.3376	&	3.0315	&	2.8117	&	2.6237	\\
&	70	&	3.3376	&	3.0377	&	2.8363	&	2.6833	\\
&	80	&	3.3376	&	3.0422	&	2.8571	&	2.7337	\\
&	90	&	3.3376	&	3.0452	&	2.8657	&	2.7518	\\
\hline
\end{tabular}
\label{table:SS2ah10k01}
\end{table}

\begin{table}[htpb]
\renewcommand\arraystretch{1.2}
\caption{Fundamental frequency $\omega b^2/h \sqrt{\rho_c/E_c}$ for simply supported Si$_3$N$_4$/SUS304 FGM square plate. $^\dag$denotes change in trend.} \centering
\begin{tabular}{cccccc}
\hline
gradient& Crack & \multicolumn{4}{c}{Crack length, $d/a$.} \\
\cline{3-6}
index, $n$ & orientation, $\theta$ & 0 & 0.4 & 0.6 & 0.8 \\
\hline
\multirow{11}{*}{2}&	0	&	3.0016	&	2.7383	&	2.5769	&	2.4747	\\
&	10	&	3.0016	&	2.7356	&	2.5692	&	2.4583	\\
&	20	&	3.0016	&	2.7315	&	2.5504	&	2.4130	\\
&	30	&	3.0016	&	2.7259	&	2.5283	&	2.3594	\\
&	40	&	3.0016	&	2.7220	&	2.5136	&	2.3223	\\
&{\bf 45}$^\dag$ &{\bf 3.0016} &{\bf 2.7202} &{\bf 2.5120} &{\bf 2.3170}	\\
&	50	&	3.0016	&	2.7219	&	2.5135	&	2.3222	\\
&	60	&	3.0016	&	2.7259	&	2.5283	&	2.3594	\\
&	70	&	3.0016	&	2.7315	&	2.5504	&	2.4130	\\
&	80	&	3.0016	&	2.7356	&	2.5692	&	2.4583	\\
&	90	&	3.0016	&	2.7383	&	2.5770	&	2.4747	\\
\cline{2-6}
\multirow{11}{*}{5}& 0	& 2.7221 &	2.4833 & 2.3371	& 2.2445	\\
& 10 & 2.7221 &	2.4809 & 2.3302	& 2.2297	\\
& 20 & 2.7221 &	2.4772 & 2.3131	& 2.1887	\\
& 30 & 2.7221 &	2.4722 & 2.2932	& 2.1402	\\
& 40 & 2.7221 &	2.4686 & 2.2798	& 2.1067	\\
&{\bf 45}$^\dag$ & {\bf 2.7221} &{\bf 2.4670} & {\bf 2.2785} &{\bf 2.1019}	\\
& 50 &	2.7221	& 2.4685 & 2.2798 &	2.1066	\\
& 60 &	2.7221	& 2.4722 & 2.2932 &	2.1402	\\
& 70 &	2.7221	& 2.4772 & 2.3132 &	2.1887	\\
& 80 &	2.7221	& 2.4809 & 2.3301 &	2.2297	\\
& 90 &	2.7221	& 2.4833 & 2.3371 &	2.2445	\\
\hline
\end{tabular}
\label{table:SS2ah10k25}
\end{table}

\subsection*{Effect of crack location}
\begin{figure}[htpb]
\centering
\subfigure[]{\includegraphics[scale=0.6]{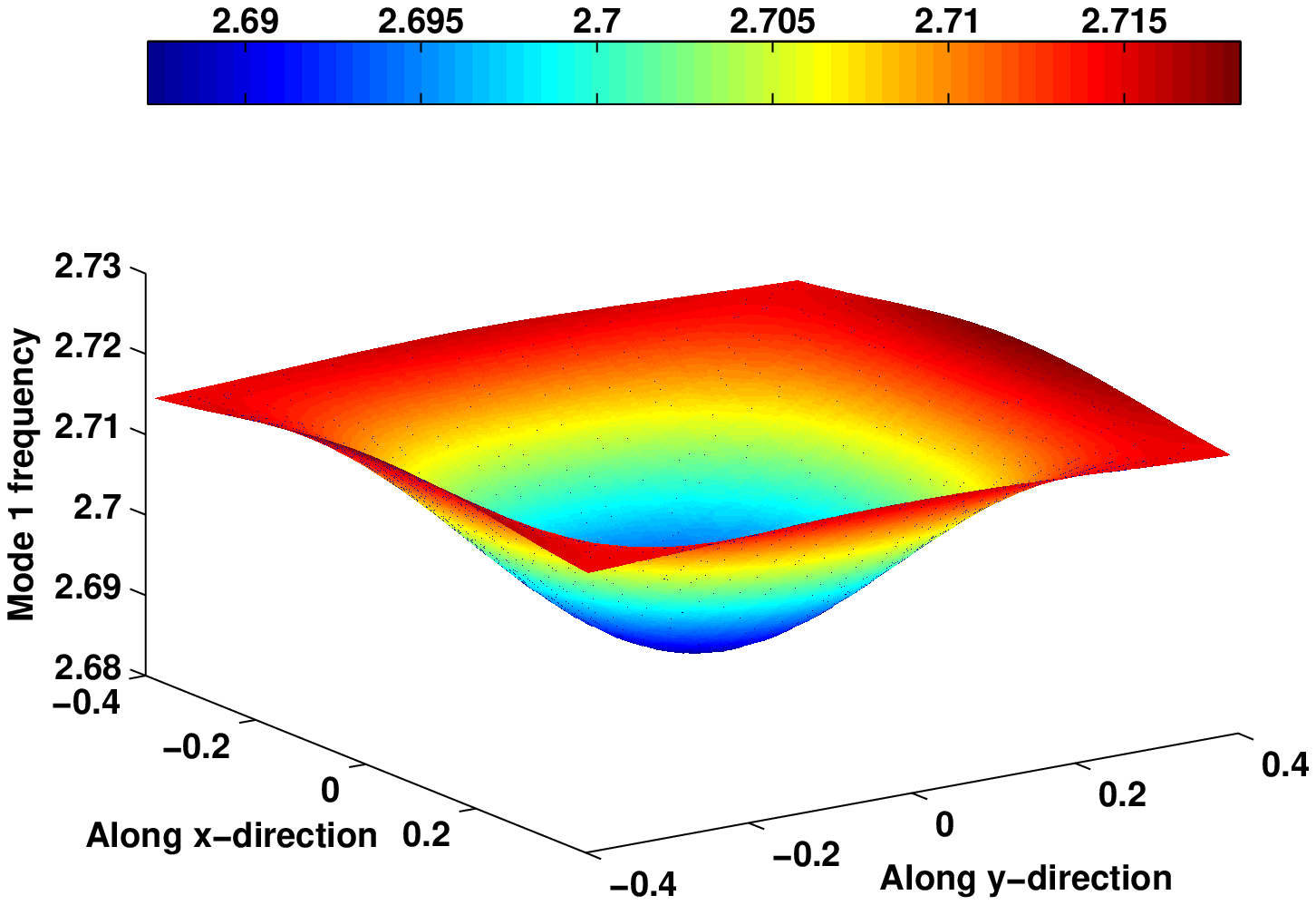}}
\subfigure[]{\includegraphics[scale=0.6]{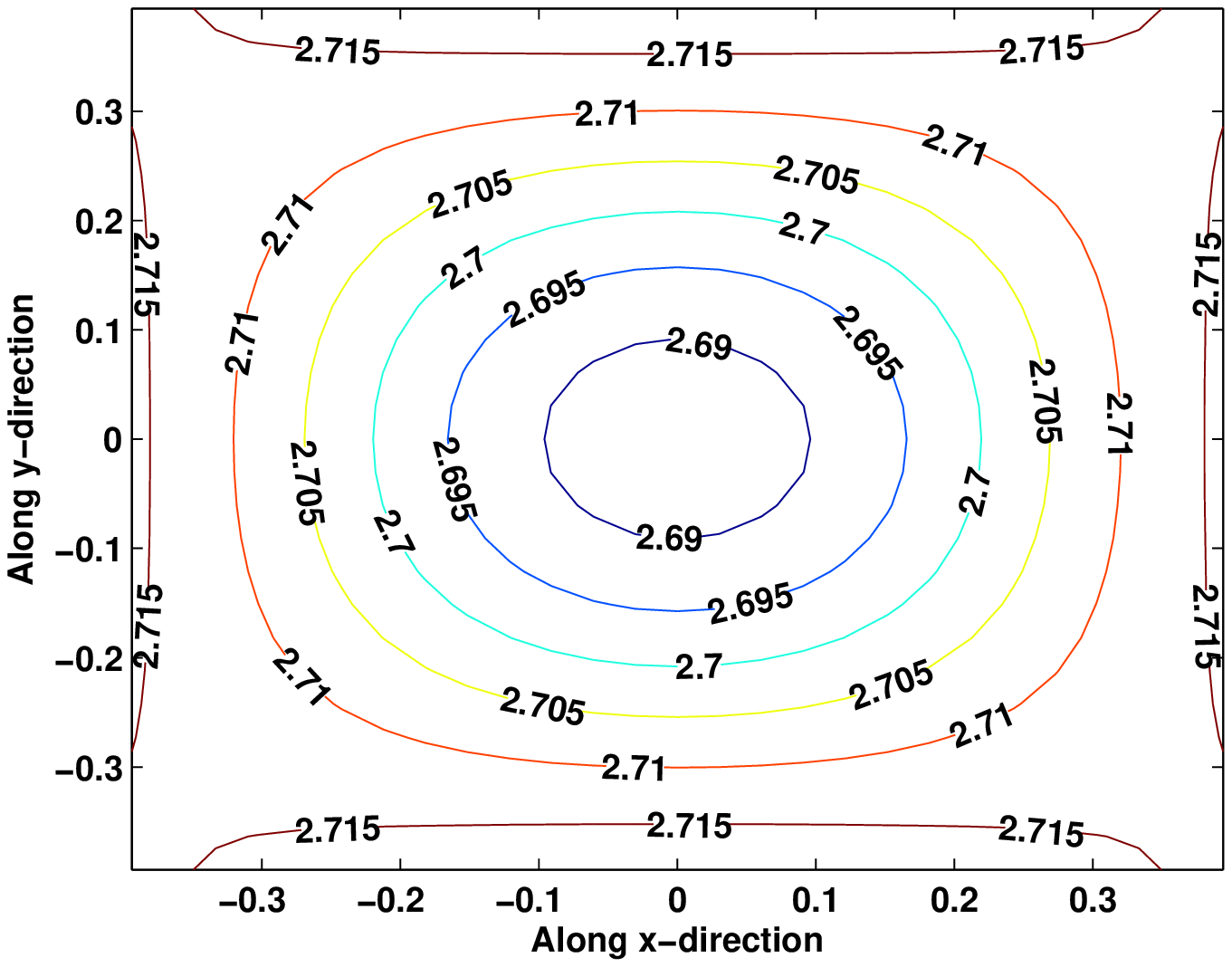}}
\caption{Variation of fundamental frequency as a function of crack position for a simply supported square FGM plate. The crack orientation, $\theta$ is taken to be 0, i.e., horizontal crack.}
\label{fig:SS2_crkLoc}
\end{figure}

Next, the influence of the crack location on the natural frequency of a square plate with thickness, $a/h=10$ and crack length, $d/a=0.2$ is studied. The results are presented in \fref{fig:SS2_crkLoc}. It is observed that the natural frequency of the plate monotonically decreases as the crack moves along the edges and towards the center of these edges. The natural frequency of the plate is maximum when the damage is situated at the corner. As the crack moves along the center lines of the plate from the edges and towards the center of the plate, the natural frequency increases up to a certain distance and then decreases. When the crack is situated at the center of the plate, the frequency is minimum.

%% file: simplyS.pstex_t
\begin{picture}(0,0)%
\includegraphics{simplyS.pstex}%
\end{picture}%
\setlength{\unitlength}{4144sp}%
\begingroup\makeatletter\ifx\SetFigFont\undefined%
\gdef\SetFigFont#1#2#3#4#5{%
  \reset@font\fontsize{#1}{#2pt}%
  \fontfamily{#3}\fontseries{#4}\fontshape{#5}%
  \selectfont}%
\fi\endgroup%
\begin{picture}(2810,2304)(2839,-4153)
\put(3331,-4066){\makebox(0,0)[b]{\smash{{\SetFigFont{12}{14.4}{\familydefault}{\mddefault}{\updefault}$x$}}}}
\put(4276,-3976){\makebox(0,0)[lb]{\smash{{\SetFigFont{12}{14.4}{\familydefault}{\mddefault}{\updefault}$a$}}}}
\put(5581,-2806){\rotatebox{90.0}{\makebox(0,0)[lb]{\smash{{\SetFigFont{12}{14.4}{\familydefault}{\mddefault}{\updefault}$b$}}}}}
\put(4321,-2671){\makebox(0,0)[lb]{\smash{{\SetFigFont{12}{14.4}{\familydefault}{\mddefault}{\updefault}$d$}}}}
\put(2926,-3796){\makebox(0,0)[lb]{\smash{{\SetFigFont{12}{14.4}{\familydefault}{\mddefault}{\updefault}$y$}}}}
\put(4231,-3211){\makebox(0,0)[lb]{\smash{{\SetFigFont{12}{14.4}{\familydefault}{\mddefault}{\updefault}$c_y$}}}}
\end{picture}%

%% file: SS2_ah10_multicrk.tex
\subsection*{Plate with multiple cracks}

\fref{fig:SS2MultiCrk} shows a plate with two cracks with lengths $a_1$ and $a_2$ and with orientations $\theta_1$ and $\theta_2$ they subtend with the horizontal. The effect of respective crack orientation on the fundamental frequency for a simply supported FGM plate is numerically studied. The horizontal and vertical separation (see \fref{fig:SS2MultiCrk}) between the crack tips is set to a constant value, $H=$0.2 and $V=$0.1, respectively. Table~\ref{table:SS2ah10k5} shows the variation of fundamental frequency for plate with gradient index $n=$5 and crack lengths $a_1=a_2=$ 0.2 as a function of crack orientations. \fref{fig:SS2multi} shows the variation of frequency as a function of orientation of one crack for different orientations of the second crack. It can be seen that with increase in crack orientation, the frequency initially decreases until it reaches a mimimum at $\theta = 45^o$. With further increase in crack orientation, the frequency increases and reaches maximum at $\theta_1=\theta_2=$90$^o$. The frequency value at $\theta =0^o$ is not the same as that at $\theta=90^o$ as we would expect. This is because when $\theta_1 = \theta_2=90^o$, the crack is located away from the center of the plate and crack disturbs the mode shape slightly.

\begin{figure}[htpb]
\centering
\input{SS2_multicrk.pstex_t}
\caption{Plate with multiple cracks: geometry}
\label{fig:SS2MultiCrk}
\end{figure}

\begin{landscape}
\begin{table}[htpb]
\renewcommand\arraystretch{1.5}
\caption{Fundamental frequency $\omega b^2/h \sqrt{\rho_c/E_c}$ for simply supported Si$_3$N$_4$/SUS304 FGM square plate, gradient index, $n=$5. $^\dag$denotes change in trend.} \centering
\begin{tabular}{cccccccccccc}
\hline
Crack & \multicolumn{11}{c}{Crack orientation, $\theta_2$.} \\
\cline{2-12}
Orientation, $\theta_1$& 0 & 10 & 20 & 30 & 40 & 45 & {\bf 50}$^\dag$& 60 & 70 & 80 & 90 \\
\hline
0 & 2.6149	&2.6122	&2.6077	&2.6034	&2.5997	&2.5988	&{\bf 2.5981}&2.5993	&2.6023	&2.6066	&2.6098\\
10	&2.6122	&2.6094	&2.6049	&2.6006	&2.5969	&2.5960	&{\bf 2.5953}&2.5964	&2.5994	&2.6037	&2.6068\\
20	&2.6077	&2.6049	&2.6004	&2.5961	&2.5924	&2.5914	&{\bf 2.5907}&2.5919	&2.5948	&2.5991	&2.6022\\
30	&2.6034	&2.6006	&2.5961	&2.5918	&2.5881	&2.5871	&{\bf 2.5864}&2.5876	&2.5905	&2.5948	&2.5979\\
40	&2.5997	&2.5969	&2.5924	&2.5881	&2.5844	&2.5835	&{\bf 2.5827}&2.5839	&2.5868	&2.5910	&2.5941\\
45	&2.5987	&2.5959	&2.5914	&2.5871	&2.5834	&2.5825	&{\bf 2.5818}&2.5829	&2.5858	&2.5901	&2.5931\\
{\bf 50}$^\dag$&{\bf 2.5980} &{\bf 2.5952} &{\bf 2.5907} & {\bf 2.5864 }& {\bf 2.5827} & {\bf 2.5818}	&{\bf 2.5811}& {\bf 2.5822} & {\bf 2.5851} &{\bf 2.5893} & {\bf 2.5924}\\
60	&2.5992	&2.5964	&2.5918	&2.5876	&2.5839	&2.5829	&{\bf 2.5822}&2.5834	&2.5863	&2.5905	&2.5935\\
70	&2.6022	&2.5993	&2.5948	&2.5905	&2.5868	&2.5858	&{\bf 2.5851}&2.5863	&2.5891	&2.5933	&2.5963\\
80	&2.6065	&2.6036	&2.5990	&2.5947	&2.5910	&2.5900	&{\bf 2.5893}&2.5904	&2.5933	&2.5974	&2.6004\\
90	&2.6098	&2.6069	&2.6022	&2.5979	&2.5942	&2.5932	&{\bf 2.5924}&2.5935	&2.5964	&2.6005	&2.6035\\
\hline
\end{tabular}
\label{table:SS2ah10k5}
\end{table}

\end{landscape}

\begin{figure}[htpb]
\centering
\subfigure[Effect of crack orientation, $\theta$]{\includegraphics[scale=0.5]{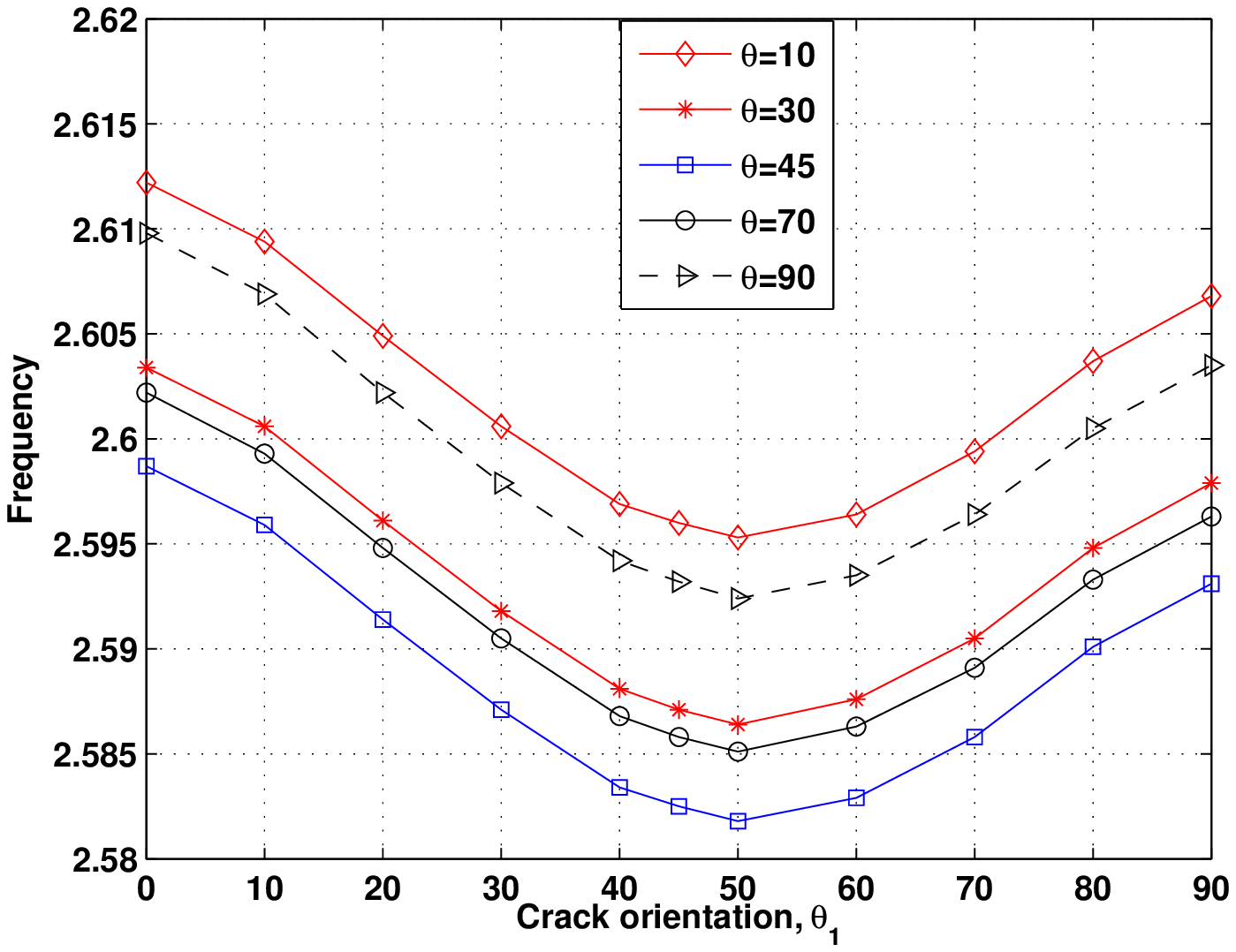}}
\subfigure[Effect of gradient index, $n$]{\includegraphics[scale=0.5]{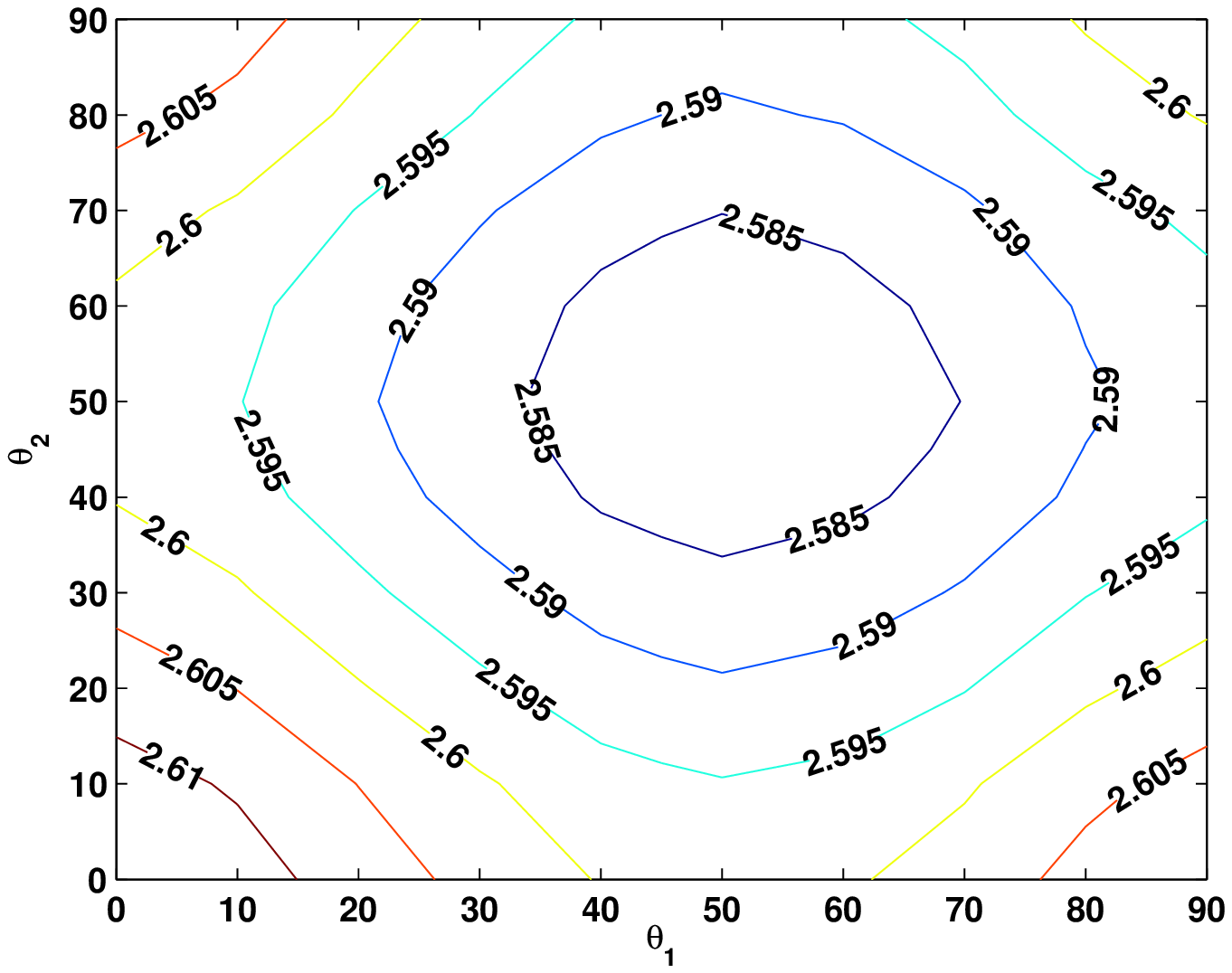}}
\caption{Variation of fundamental frequency for a simply supported square Si$_3$N$_4$/SUS304 FGM plate with a center crack as a function crack orientation with $H=$0.2, $V=$0.1 and $a/h=$10 (see \fref{fig:SS2MultiCrk}).}
\label{fig:SS2multi}
\end{figure}

%% file: SS2_multicrk.pstex_t
\begin{picture}(0,0)%
\includegraphics{SS2_multicrk.pstex}%
\end{picture}%
\setlength{\unitlength}{2279sp}%
\begingroup\makeatletter\ifx\SetFigFont\undefined%
\gdef\SetFigFont#1#2#3#4#5{%
  \reset@font\fontsize{#1}{#2pt}%
  \fontfamily{#3}\fontseries{#4}\fontshape{#5}%
  \selectfont}%
\fi\endgroup%
\begin{picture}(8947,7736)(1291,-7190)
\put(3421,-2581){\makebox(0,0)[lb]{\smash{{\SetFigFont{12}{14.4}{\familydefault}{\mddefault}{\updefault}$\theta_2$}}}}
\put(6211,-4561){\makebox(0,0)[lb]{\smash{{\SetFigFont{12}{14.4}{\familydefault}{\mddefault}{\updefault}$a_1$}}}}
\put(6526,-4021){\makebox(0,0)[lb]{\smash{{\SetFigFont{12}{14.4}{\familydefault}{\mddefault}{\updefault}$\theta_1$}}}}
\put(4996,-7081){\makebox(0,0)[lb]{\smash{{\SetFigFont{12}{14.4}{\familydefault}{\mddefault}{\updefault}$a$}}}}
\put(1306,-3616){\makebox(0,0)[lb]{\smash{{\SetFigFont{12}{14.4}{\familydefault}{\mddefault}{\updefault}$b$}}}}
\put(9766,-6226){\makebox(0,0)[lb]{\smash{{\SetFigFont{12}{14.4}{\familydefault}{\mddefault}{\updefault}$x$}}}}
\put(2341,164){\makebox(0,0)[lb]{\smash{{\SetFigFont{12}{14.4}{\familydefault}{\mddefault}{\updefault}$y$}}}}
\put(3691,-1996){\makebox(0,0)[lb]{\smash{{\SetFigFont{12}{14.4}{\familydefault}{\mddefault}{\updefault}$a_2$}}}}
\put(5986,-3301){\makebox(0,0)[lb]{\smash{{\SetFigFont{12}{14.4}{\familydefault}{\mddefault}{\updefault}$V$}}}}
\put(4951,-4696){\makebox(0,0)[lb]{\smash{{\SetFigFont{12}{14.4}{\familydefault}{\mddefault}{\updefault}$H$}}}}
\end{picture}%

%% file: cantilever.pstex_t
\begin{picture}(0,0)%
\includegraphics{cantilever.pstex}%
\end{picture}%
\setlength{\unitlength}{2279sp}%
\begingroup\makeatletter\ifx\SetFigFont\undefined%
\gdef\SetFigFont#1#2#3#4#5{%
  \reset@font\fontsize{#1}{#2pt}%
  \fontfamily{#3}\fontseries{#4}\fontshape{#5}%
  \selectfont}%
\fi\endgroup%
\begin{picture}(7291,5497)(2187,-5885)
\put(2202,-3121){\makebox(0,0)[lb]{\smash{{\SetFigFont{12}{14.4}{\familydefault}{\mddefault}{\updefault}$b$}}}}
\put(5401,-5776){\makebox(0,0)[lb]{\smash{{\SetFigFont{12}{14.4}{\familydefault}{\mddefault}{\updefault}$a$}}}}
\put(8101,-4111){\makebox(0,0)[lb]{\smash{{\SetFigFont{12}{14.4}{\familydefault}{\mddefault}{\updefault}$c_y$}}}}
\put(6104,-2643){\makebox(0,0)[lb]{\smash{{\SetFigFont{12}{14.4}{\familydefault}{\mddefault}{\updefault}$d$}}}}
\put(6329,-3150){\makebox(0,0)[lb]{\smash{{\SetFigFont{12}{14.4}{\familydefault}{\mddefault}{\updefault}$\theta$}}}}
\put(9070,-4908){\makebox(0,0)[lb]{\smash{{\SetFigFont{12}{14.4}{\familydefault}{\mddefault}{\updefault}$x$}}}}
\put(3225,-691){\makebox(0,0)[lb]{\smash{{\SetFigFont{12}{14.4}{\familydefault}{\mddefault}{\updefault}$y$}}}}
\end{picture}%

%% file: 2011-03_FGM_vib.bbl
\begin{thebibliography}{10}
\expandafter\ifx\csname url\endcsname\relax
  \def\url#1{\texttt{#1}}\fi
\expandafter\ifx\csname urlprefix\endcsname\relax\def\urlprefix{URL }\fi
\expandafter\ifx\csname href\endcsname\relax
  \def\href#1#2{#2} \def\path#1{#1}\fi

\bibitem{Koizumi1993}
M.~Koizumi, The concept of {FGM}, Ceram. Trans. Funct. Graded Mater 34 (1993)
  3--10.

\bibitem{Koizumi1997}
M.~Koizumi, {FGM} activities in {J}apan, Composites 28 (1997) 1--4.

\bibitem{Vel2004}
S.~Vel, R.~Batra, Three-dimensional exact solution for the vibration of
  functionally graded rectangular plates, Journal of Sound and Vibration 272
  (2004) 703--730.

\bibitem{Matsunaga2008}
H.~Matsunaga, Free vibration and stability of functionally graded plates
  according to a 2-{D} higher-order deformation theory, Composite Structures 82
  (2008) 499--512.

\bibitem{He2001}
X.~He, T.~Ng, S.~Sivashankar, K.~Liew, Active control of fgm plates with
  integrated piezoelectric sensors and actuators, International Journal of
  Solids and Structures 38 (2001) 1641--1655.

\bibitem{Liew2001}
K.~Liew, X.~He, T.~Ng, S.~Sivashankar, Active control of fgm plates subjected
  to a temperature gradient: modeling via finite element method based on
  {FSDT}, International Journal for Numerical Methods in Engineering 52 (2001)
  1253--1271.

\bibitem{Ng2000}
T.~Ng, K.~Lam, K.~Liew, Effect of {FGM} materials on parametric response of
  plate structures, Computer Methods in Applied Mechanics and Engineering 190
  (2000) 953--962.

\bibitem{Yang2002}
J.~Yang, H.~Shen, Vibration characteristic and transient response of
  shear-deformable functionally graded plates in thermal environments, Journal
  of Sound and Vibration 255 (2002) 579--602.

\bibitem{Ferreira2006}
A.~Ferreira, R.~Batra, C.~Roque, L.~Qian, R.~Jorge, Natural frequencies of
  functionally graded plates by a meshless method, Composite Structures 75
  (2006) 593--600.

\bibitem{Yang2001}
J.~Yang, H.~Shen, Dynamic response of initially stressed functionally graded
  rectangular thin plates, Composite Structures 54 (2001) 497--508.

\bibitem{Qian2004}
L.~Qian, R.~Batra, L.~Chen, Static and dynamic deformations of thick
  functionally graded elastic plates by using higher order shear and normal
  deformable plate theory and meshless local petrov galerkin method, Composites
  Part B: Engineering 35 (2004) 685--697.

\bibitem{Reddy2000}
J.~Reddy, Analysis of functionally graded plates, International Journal for
  Numerical Methods in Engineering 47 (2000) 663--684.

\bibitem{R2010}
A.~A. R, A.~Bagri, S.~Bordas, T.~Rabczuk, Analysis of thermoelastic waves in a
  two-dimensional functionally graded materials domain by the meshless local
  {P}etrov-{G}alerkin method, Computer Modelling in Engineering and Science 65
  (2010) 27--74.

\bibitem{Lynn1967}
P.~Lynn, N.~Kumbasar, Free vibrations of thin rectangular plates having narrow
  cracks with simply supported edges, Dev. Mech 4 (1967) 928--991.

\bibitem{Stahl1972}
B.~Stahl, L.~Keer, Vibration and stability of cracked rectangular plates.,
  International Journal of Solids and Structures 8 (1972) 69--91.

\bibitem{Solecki1985}
R.~Solecki, Bending vibration of rectangular plate with arbitrarly located
  rectilinear crack, Engineering Fracture Mechanics 22~(4) (1985) 687--695.

\bibitem{Khadem2000}
S.~Khadem, M.~Rezaee, Introduction of modified comparison functions for
  vibration analysis of a rectangular cracked plate., Journal of Sound and
  Vibration 236~(2) (2000) 245--258.

\bibitem{Wu2005}
G.~Wu, Y.~Shih, Dynamic instability of rectangular plate with an edge cracked,
  Computers and Structures 84~(1-2) (2005) 1--10.

\bibitem{Qian1991}
G.~Qian, S.~Gu, J.~Jiang, A finite element model of cracked plates and
  application to vibration problems, Computers and Structures 39~(5) (1991)
  483--487.

\bibitem{Lee1993}
H.~Lee, S.~Lim, Vibration of cracked rectangular plates including transverse
  shear deformation and rotary inertia, Computers and Structures 49~(4) (1993)
  715--718.

\bibitem{Huang2011}
C.~Huang, O.~M. III, M.~Chang, Vibrations of cracked rectangular {FGM} thick
  plates, Commposite Structures\href
  {http://dx.doi.org/10.1016/j.compstruct.2011.01.005}
  {\path{doi:10.1016/j.compstruct.2011.01.005}}.

\bibitem{Kitipornchai2009}
S.~Kitipornchai, L.~Ke, J.~Y. andY Xiang, Nonlinear vibration of edge cracked
  functionally graded {T}imoshenko beams, Journal of Sound and Vibration 324
  (2009) 962--982.

\bibitem{Yang2010}
J.~Yang, Y.~Hao, W.~Zhang, S.~Kitipornchai, Nonlinear dynamic response of a
  functionally graded plate with a through-width surface crack, Nonlinear
  {D}ynamics 59 (2010) 207--219.

\bibitem{Somashekar1987}
B.~Somashekar, G.~Prathap, C.~R. Babu, A field-consistent four-noded laminated
  anisotropic plate/shell element, Computers and Structures 25 (1987) 345--353.

\bibitem{Bachene2009}
M.~Bachene, R.~Tiberkak, S.~Rechak, Vibration analysis of cracked plates using
  the extended finite element method, Arch. Appl. Mech. 79 (2009) 249--262.
\newblock \href {http://dx.doi.org/10.1007/s00419-008-0224-7}
  {\path{doi:10.1007/s00419-008-0224-7}}.

\bibitem{Bachene2009a}
M.~Bachene, R.~Tiberkak, S.~Rechak, G.~Maurice, B.~Hachi, Enriched finite
  element for modal analysis of cracked plates, in: Damage and {F}racture
  {M}echanics: {F}ailure {A}nalysis and {E}ngineering {M}aterials and
  {S}tructures, Springer Science, 2009, pp. 463--471.

\bibitem{Tiberkak2009}
R.~Tiberkak, M.~Bachene, B.~Hachi, S.~Rechak, M.~Haboussi, Dynamic response of
  cracked plate subjected to impact loading using the extended finite element
  method, in: Damage and {F}racture {M}echanics: {F}ailure {A}nalysis and
  {E}ngineering {M}aterials and {S}tructures, Springer Science, 2009, pp.
  297--306.

\bibitem{Mori1973}
T.~Mori, K.~Tanaka, Average stress in matrix and average elasticenergy of
  materials with misfitting inclusions, Acta Metallurgica 21 (1973) 571--574.

\bibitem{Benvensite1987}
Y.~Benvensite, A new approach to the application of mori–tanaka's theory in
  composite materials, Mechanics of Materials 6 (1987) 147--157.

\bibitem{Cheng2000}
Z.-Q. Cheng, R.~Batra, Three dimensional thermoelastic deformations of a
  functionally graded elliptic plate, Composites Part B: Engineering 2 (2000)
  97--106.

\bibitem{Reddy1998}
J.~Reddy, C.~Chin, Thermomechanical analysis of functionally graded cylinders
  and plates, Journal of Thermal Stresses 21 (1998) 593--629.

\bibitem{Rajasekaran1973}
S.~Rajasekaran, D.~Murray, Incremental finite element matrices, {ASCE} Journal
  of Structural Divison 99 (1973) 2423--2438.

\bibitem{Ganapathi1991}
M.~Ganapathi, T.~Varadan, B.~Sarma, Nonlinear flexural vibrations of laminated
  orthotropic plates, Computers and Structures 39 (1991) 685--688.

\bibitem{Babuvska1994}
I.~Babu\v{s}ka, G.~Caloz, J.~Osborn, Special finite element methods for a class
  of second order elliptic problems with rough coefficients, SIAM Journal of
  Numerical Analysis 31 (1994) 945--981.

\bibitem{Belytschko1999}
T.~Belytschko, T.~Black, Elastic crack growth in finite elements with minimal
  remeshing, International Journal for Numerical Methods in Engineering 45
  (1999) 601--620.

\bibitem{Melenk1995}
J.~M. Melenk, On generalized finite element methods, Ph.D. thesis, University
  of Maryland, College Park, MD (1995).

\bibitem{Duarte2007}
C.~A. Duarte, T.~J. Liszka, W.~W. Tworzydlo, Clustered generalized finite
  element methods for mesh unrefinement, non-matching and invalid meshes,
  International Journal for Numerical Methods in Engineering 69~(11) (2007)
  2409--2440.
\newblock \href {http://dx.doi.org/10.1002/nme.1862}
  {\path{doi:10.1002/nme.1862}}.

\bibitem{Babuvska1997}
I.~Babu\v{s}ka, J.~Melenk, The partition of unity finite element method,
  International Journal for Numerical Methods in Engineering 40 (1997)
  727--758.

\bibitem{Babuska2008}
I.~Babuska, V.~Nistor, N.~Tarfulea, Generalized finite element method for
  second-order elliptic operators with dirichlet boundary conditions, Journal
  of Computational and Applied Mathematics 218~(1) (2008) 175--183.
\newblock \href {http://dx.doi.org/10.1016/j.cam.2007.04.041}
  {\path{doi:10.1016/j.cam.2007.04.041}}.

\bibitem{Simone2006}
A.~Simone, C.~Duarte, E.~V. der Giessen, A generalized finite element method
  for polycrystals with discontinuous grain boundaries, International Journal
  for Numerical Methods in Engineering 67 (2006) 1122--1145.

\bibitem{Karihaloo2003}
B.~Karihaloo, Q.~Xiao, Modelling of stationary and growing cracks in {FE}
  framework without remeshing: a state-of-the-art review, Computers and
  Structures 81~(3) (2003) 119--129.
\newblock \href {http://dx.doi.org/10.1016/S0045-7949(02)00431-5}
  {\path{doi:10.1016/S0045-7949(02)00431-5}}.

\bibitem{Yazid2009}
A.~Yazid, N.~Abdelkader, H.~Abdelmadjid, A state-of-the-art review of the x-fem
  for computational fracture mechanics, Applied Mathematical Modelling 33~(12)
  (2009) 4269--4282.
\newblock \href {http://dx.doi.org/10.1016/j.apm.2009.02.010}
  {\path{doi:10.1016/j.apm.2009.02.010}}.

\bibitem{Bordas2007a}
S.~Bordas, P.~V. Nguyen, C.~Dunant, A.~Guidoum, H.~Nguyen-Dang, An extended
  finite element library, International Journal for Numerical Methods in
  Engineering 71 (2007) 703--732.

\bibitem{Rabczuk2010}
T.~Rabczuk, S.~Bordas, G.~Zi, On three-dimensional modelling of crack growth
  using partition of unity methods, Computers and Structures 88 (2010)
  1391--1411.

\bibitem{Dolbow2000}
J.~Dolbow, N.~Mo{\"e}s, T.~Belytschko, Modeling fracture in mindlin-reissner
  plates with the extended finite element method, International Journal of
  Solids and Structures 37~(48-50) (2000) 7161--7183.
\newblock \href {http://dx.doi.org/10.1016/S0020-7683(00)00194-3}
  {\path{doi:10.1016/S0020-7683(00)00194-3}}.

\bibitem{Moes1999}
N.~Mo{\"e}s, J.~Dolbow, T.~Belytschko, A finite element method for crack growth
  without remeshing, International Journal for Numerical Methods in Engineering
  46~(1) (1999) 131--150.

\bibitem{Watkins2002}
D.~S. Watkins, Fundamentals of matrix computations, Wiley Publishers, 2002.

\bibitem{Sundararajan2005}
N.~Sundararajan, T.~Prakash, M.~Ganapathi, Nonlinear free flexural vibrations
  of functionally graded rectangular and skew plates under thermal
  environments, Finite Elements in Analysis and Design 42~(2) (2005) 152--168.

\bibitem{Prakash2007}
T.~Prakash, N.~Sundararajan, M.~Ganapathi, On the nonlinear axisymmetric
  dynamic buckling behavior of clamped functionally graded spherical caps,
  Journal of Sound and Vibration 299 (2007) 36--43.

\bibitem{Singh2011}
M.~Singh, T.~Prakash, M.~Ganapathi, Finite element analysis of functionally
  graded plates under transverse load, Finite Elements in Analysis and Design
  47 (2011) 453--460.

\bibitem{Huang2009}
C.~Huang, A.~Lieissa, Vibration analysis of rectangular plates with side cracks
  via the {R}itz method, Journal of Sound and Vibration 323 (2009) 974--988.

\end{thebibliography}
